\documentclass[leqno]{article}

\usepackage[frenchb,english]{babel}
\usepackage[T1]{fontenc}

\usepackage{amsmath}
\usepackage{amssymb}
\usepackage{amsfonts}
\usepackage{enumerate}
\usepackage{vmargin}
\usepackage[all]{xy}
\usepackage{mathrsfs}
\usepackage{mathtools}
\setmarginsrb{3cm}{4cm}{2.5cm}{2cm}{0cm}{0cm}{1.5cm}{3cm}
\newcommand{\tr}{\ensuremath{\mathop{\rm Tr\,}\nolimits}}

\renewcommand{\leq}{\ensuremath{\leqslant}}
\renewcommand{\geq}{\ensuremath{\geqslant}}
\newcommand{\n}{\noindent}
\newcommand{\qed}{\hfill \vrule height6pt  width6pt depth0pt}

\newcommand\be{\begin{eqnarray}}
\newcommand\ee{\end{eqnarray}}
\newcommand{\norm}[1]{ \| #1  \|}
\newcommand{\bnorm}[1]{ \big\| #1  \big\|}
\newcommand{\Bnorm}[1]{ \Big\| #1  \Big\|}
\newcommand{\bgnorm}[1]{ \bigg\| #1  \bigg\|}
\newcommand{\Bgnorm}[1]{ \Bigg\| #1  \Bigg\|}
\newcommand{\xra}{\xrightarrow}

\newcommand{\ot}{\otimes}

\def\Ker{{\rm Ker} \, }

\def\Re{{\rm Re} }

\def\diag{{\rm diag} }

\newcommand{\ovl}{\overline}

\newcommand{\Rad}{{\rm Rad}}
\newcommand{\rad}{{\rm rad}}
\newcommand{\col}{{\rm Col}}

\newcommand{\Ran}{{\rm Ran}}

\newtheorem{thm}{Theorem}[section]
\newtheorem{defi}[thm]{Definition}
\newtheorem{prop}[thm]{Proposition}

\newtheorem{cor}[thm]{Corollary}
\newtheorem{lemma}[thm]{Lemma}

\newtheorem{remark}[thm]{Remark}

\newenvironment{preuve}[1][]{\noindent {\it Proof #1} : }{\hbox{~}\qed
\smallskip
}

\newcommand{\beq}{\begin{equation}}
\newcommand{\eeq}{\end{equation}}

\numberwithin{equation}{section}

\begin{document}
\selectlanguage{english}
\title{\bfseries{Square functions for Ritt operators on noncommutative $L^p$-spaces}}
\date{}
\author{\bfseries{C\'edric Arhancet}}

\maketitle

\begin{abstract}
For any Ritt operator $T$ acting on a noncommutative $L^p$-space, we
define the notion of \textit{completely} bounded functional calculus
$H^\infty(B_\gamma)$ where $B_\gamma$ is a Stolz domain. Moreover,
we introduce the `column square functions'
$\norm{x}_{p,T,c,\alpha}=\Bnorm{\Big(\sum_{k=1}^{+\infty}k^{2\alpha-1}\left|T^{k-1}(I-T)^{\alpha}(x)\right|^2\Big)^{\frac{1}{2}}}_{L^p(M)}$
and the `row square functions'
$\norm{x}_{p,T,r,\alpha}=\Bnorm{\Big(\sum_{k=1}^{+\infty}k^{2\alpha-1}
\left|\Big(T^{k-1}(I-T)^{\alpha}(x)\Big)^*\right|^2\Big)^{\frac{1}{2}}}_{L^p(M)}$
for any $\alpha>0$ and any $x\in L^p(M)$. Then, we provide an
example of Ritt operator which admits a completely bounded
$H^\infty(B_\gamma)$ functional calculus for some $\gamma \in
\big]0,\frac{\pi}{2}\big[$ such that the square functions
$\norm{\cdot}_{p,T,c,\alpha}$ and $\norm{\cdot}_{p,T,r,\alpha}$ are
not equivalent. Moreover, assuming $1<p<2$ and $\alpha>0$, we prove
that if $\Ran (I-T)$ is dense and $T$ admits a completely bounded
$H^\infty(B_\gamma)$ functional calculus for some $\gamma \in
\big]0,\frac{\pi}{2}\big[$ then there exists a positive constant $C$
such that for any $x \in L^p(M)$, there exists $x_1, x_2 \in L^p(M)$
satisfying $x=x_1+x_2$ and
$\norm{x_1}_{p,T,c,\alpha}+\norm{x_2}_{p,T,r,\alpha}\leq C
\norm{x}_{L^p(M)}$. Finally, we observe that this result applies to
a suitable class of selfadjoint Markov maps on noncommutative
$L^p$-spaces.
\bigskip

\end{abstract}


\makeatletter
 \renewcommand{\@makefntext}[1]{#1}
 \makeatother \footnotetext{\noindent
 This work is partially supported by ANR 06-BLAN-0015.
}


\section{Introduction}


\indent Let $M$ be a semifinite von Neumann algebra equipped with a
normal semifinite faithful trace. For any $1\leq p <\infty$, we let
$L^p(M)$ denote the associated (noncommutative) $L^p$-space. 
Let $T$ be a bounded operator on $L^p(M)$. Consider the following
`square function'
\begin{equation}\label{norm T1}
\norm{x}_{p,T,1}=\inf\Bigg\{\Bgnorm{\bigg(\sum_{k=1}^{+\infty}
|u_k|^2\bigg)^{\frac{1}{2}}}_{L^p}+\Bgnorm{\bigg(\sum_{k=1}^{+\infty}|v_k^*|^2\bigg)^{\frac{1}{2}}}_{L^p}
\colon\ u_k+v_k=k^{\frac{1}{2}}\big(T^{k}(x)-T^{k-1}(x)\big)\text{
for any $k$}\Bigg\}
\end{equation}
if $1<p\leq 2$ and
\begin{equation}\label{norm T1 prime}
\norm{x}_{p,T,1}=\max \Bigg\{
\Bgnorm{\bigg(\sum_{k=1}^{+\infty}k\left|T^{k}(x)-T^{k-1}(x)\right|^2\bigg)^{\frac{1}{2}}}_{L^p}
,\Bgnorm{\bigg(\sum_{k=1}^{+\infty}k\left|\Big(T^{k}(x)-T^{k-1}(x)\Big)^*\right|^2\bigg)^{\frac{1}{2}}}_{L^p}\Bigg\}
\end{equation}
if $2\leq p<\infty$, defined for any $x \in L^p(M)$. Such quantities
were introduced in \cite{LM0} and studied in this paper and in
\cite{ALM}. Similar square functions for continuous semigroups
played a key role in the recent development of $H^\infty$-calculus
and its applications. See in particular the paper \cite{JMX}, the
survey \cite{LM3} and the references therein.

For any $\gamma\in \big]0,\frac{\pi}{2}\big[$, let $B_\gamma$ be the
interior of the convex hull of $1$ and the disc $D(0,\sin \gamma)$.
Suppose $1 < p <\infty$. Let $T$ be a Ritt operator with $\Ran(I-T)$
dense in $L^p(M)$ which admits a bounded $H^\infty(B_\gamma)$
functional calculus for some $\gamma\in \big]0,\frac{\pi}{2}\big[$,
i.e. there exists an angle $\gamma\in \big]0,\frac{\pi}{2}\big[$ and
a positive constant $K$ such that $\bnorm{\varphi(T)}_{L^p(M) \to
L^p(M)} \leq K \norm{\varphi}_{ H^\infty(B_\gamma)}$ for any complex
polynomial $\varphi$. A result of \cite{LM0} essentially says that 
\begin{equation}
\label{equivalence usual norm and T1}
 \norm{x}_{L^p(M)}\approx
\norm{x}_{p,T,1},\qquad x\in L^p(M)
\end{equation}
(see also \cite[Remark 6.4]{ALM}). Now, consider the following
`column and row square functions'
\begin{equation}
\label{column and row square functions}
\norm{x}_{p,T,c,1}=\Bgnorm{\bigg(\sum_{k=1}^{+\infty}k\left|T^{k}(x)-T^{k-1}(x)\right|^2\bigg)^{\frac{1}{2}}}_{L^p}
 \text{and}  \
\norm{x}_{p,T,r,1}=\Bgnorm{\bigg(\sum_{k=1}^{+\infty}k\left|\Big(T^{k}(x)-T^{k-1}(x)\Big)^*\right|^2\bigg)^{\frac{1}{2}}}_{L^p}
\end{equation}
defined for any $x \in L^p(M)$. Assume $1<p<2$. In this context, if
$x\in L^p(M)$, it is natural to search sufficient conditions to find
a decomposition $x=x_1+x_2$ such that $\norm{x_1}_{p,T,c,1}$ and
$\norm{x_2}_{p,T,r,1}$ are finite. The first main result of this
paper is the next theorem. It strengthens the above equivalence
(\ref{equivalence usual norm and T1}) in the case where $T$ actually
admits a \textit{completely} bounded $H^\infty(B_\gamma)$ functional
calculus, i.e. there exists a positive constant $K$ such that
$\bnorm{\varphi(T)}_{cb, L^p(M) \to L^p(M)} \leq K \norm{\varphi}_{
H^\infty(B_\gamma)}$ for any complex polynomial $\varphi$.
\begin{thm}
\label{Th 2} Suppose $1<p<2$. Let $T$ be a Ritt operator on $L^p(M)$
with $\Ran (I-T)$ dense in $L^p(M)$. Assume that $T$ admits a
completely bounded $H^\infty(B_\gamma)$ functional calculus for some
$\gamma \in \big]0,\frac{\pi}{2}\big[$. Then we have
$$
\norm{x}_{L^p(M)} \approx
\inf\Big\{\norm{x_1}_{p,T,c,1}+\norm{x_2}_{p,T,r,1}\ \colon\
x=x_1+x_2 \Big\},\qquad x \in L^p(M).
$$
\end{thm}

In this context, it is natural to compare the both quantities of
(\ref{column and row square functions}). The second principal result
of this paper is the following theorem. It says that in general,
`column and row square functions' are not equivalent.
\begin{thm}
\label{Th 3} Suppose $1<p\not=2<\infty$. Then there exists a Ritt
operator $T$ on the Schatten space $S^p$, with $\Ran(I-T)$ dense in
$S^p$, which admits a completely bounded $H^\infty(B_\gamma)$
functional calculus for some $\gamma\in \big]0,\frac{\pi}{2}\big[$
such that
\begin{equation}\label{c/r=infty}
\sup\Bigg\{\frac{\norm{x}_{p,T,c,1}}{\norm{x}_{p,T,r,1}}\ \colon\
x\in S^p\Bigg\}=\infty\  \text{if $2 < p < \infty$}\ \text{and} \
\sup\Bigg\{\frac{\norm{x}_{p,T,r,1}}{\norm{x}_{p,T,c,1}}\ \colon\
x\in S^p\Bigg\}=\infty\ \text{if $1 < p < 2$.}
\end{equation}
Moreover, the same result holds with $\norm{\cdot}_{p,T,c,1}$ and
$\norm{\cdot}_{p,T,r,1}$ switched.
\end{thm}


The paper is organized as follows. Section 2 gives a brief
presentation of noncommutative $L^p$-spaces and Ritt operators and
we introduce the notions of Col-Ritt and Row-Ritt operators and
completely bounded $H^\infty(B_\gamma)$ functional calculus which
are relevant to our paper. The next section 3 mostly contains
preliminary results concerning Col-Ritt and Row-Ritt operators.
Section 4 is devoted to prove Theorems \ref{Th 3}. In section 5, we
present a proof of Theorem \ref{Th 2}. We end this section by giving
some natural examples to which this result can be applied.

In the above presentation and later on in the paper we will use
$\lesssim$ to indicate an inequality up to a constant which does not
depend to the particular element to which it applies. Then
$A(x)\approx B(x)$ will mean that we both have $A(x)\lesssim B(x)$
and $B(x)\lesssim A(x)$.

\section{Background and preliminaries}

We start with a few preliminaries on noncommutative $L^p$-spaces.
Let $M$ be a von Neumann algebra equipped with a normal semifinite
faithful trace $\tau$. Let $M_+$ be the set of all positive elements
of $M$ and let $S_+$ be the set of all $x$ in $M_+$ such that
$\tau(x)<\infty$. Then let $S$ be the linear span of $S_+$. For any
$1\leq p<\infty$, define
$$
\norm{x}_{L^p(M)}\,=\,\bigl(\tau(\vert
x\vert^p)\bigr)^{\frac{1}{p}},\qquad x\in S,
$$
where $\vert x\vert =(x^*x)^{\frac{1}{2}}$ is the modulus of $x$.
Then $\big(S,\norm{\cdot}_{L^p(M)}\big)$ is a normed space. The
corresponding completion is the noncommutative $L^p$-space
associated with $(M,\tau)$ and is denoted by $L^p(M)$. By
convention, we set $L^\infty(M)=M$, equipped with the operator norm.
The elements of $L^p(M)$ can also be described as measurable
operators with respect to $(M,\tau)$. Further multiplication of
measurable operators leads to contractive bilinear maps
$L^{p}(M)\times L^{q}(M)\to L^{r}(M)$ for any $1 \leq
p,q,r\leq\infty$ such that $\frac{1}{p} + \frac{1}{q}= \frac{1}{r}$
(noncommutative H\"older's inequality). Using trace duality, we then
have $L^{p}(M)^*\,=\, L^{p^*}(M)$ isometrically for any $1\leq
p<\infty$. Moreover, complex interpolation yields
$L^p(M)=\big[L^\infty(M),L^1(M)\big]_{\frac{1}{p}}$ for any $1 \leq
p\leq \infty$. We refer the reader to \cite{PX} for details and
complements.

Let $1 \leq p< \infty$. If we equip the space $B(\ell^2)$ with the
operator norm and the canonical trace $\tr$, the space
$L^p\big(B(\ell^2)\big)$ identifies to the Schatten-von Neumann
class $S^p$. This is the space of those compact operators $x$ from
$\ell^2$ into $\ell^2$ such that
$\norm{x}_{S^p}=\big(\tr(x^*x)^{\frac{p}{2}}\big)^{\frac{1}{p}}<\infty$.
Elements of $B(\ell^2)$ or $S^p$ are regarded as matrices
$A=[a_{ij}]_{i,j \geq 1}$ in the usual way.

If the von Neumann algebra $B(\ell^2)\ovl{\ot}M$ is equipped with
the semifinite normal faithful trace $\tr\ot \tau$, the space
$L^p\big(B(\ell^2)\ovl{\ot}M\big)$ canonically identifies to a space
$S^p\big(L^p(M)\big)$ of matrices with entries in $L^p(M)$.
Moreover, under this identification, the algebraic tensor product
$S^p \ot L^p(M)$ is dense in $S^p\big(L^p(M)\big)$. We refer to
\cite{Pis3} for more about these spaces and complements.

If $1\leq p < \infty$, we say that a linear map on $L^p(M)$ is
completely bounded if $I_{S^p} \ot T$ extends to a bounded operator
on $S^p\big(L^p(M)\big)$. In this case, the completely bounded norm
$\norm{T}_{cb,L^p(M)\xra{}L^p(M)}$ of $T$ is defined by
$\norm{T}_{cb,L^p(M)\xra{}L^p(M)}=\bnorm{I_{S^p} \ot
T}_{S^p(L^p(M))\xra{}S^p(L^p(M))}$. We use the convention to define
$\norm{T}_{cb,L^p(M)\xra{}L^p(M)}$ by $+\infty$ if $T$ is not
completely bounded.

We shall use various $\ell^2$-valued noncommutative $L^p$ spaces. We
refer to \cite[Chapter 2]{JMX} for more information on these spaces.
For any $\sum_{k=1}^nx_k\ot a_k\in L^p(M)\ot \ell^2$, we set
$$
\Bgnorm{\sum_{k=1}^n x_k\ot
a_k}_{L^p(M,\ell^2_c)}=\Bgnorm{\bigg(\sum_{i,j=1}^{n} \langle
a_j,a_i\rangle x_i^*x_j\bigg)^{\frac{1}{2}}}_{L^p(M)}.
$$
We have for any family $(x_k)_{k\geq 1}$ in $L^{p}(M)$
\begin{equation}
\label{expression norm} \Bgnorm{\sum_{k=1}^{n} x_k\ot
e_k}_{L^p(M,\ell^2_c)}=\Bgnorm{\bigg(\sum_{k=1}^n|x_k|^2\bigg)^{\frac{1}{2}}}_{L^p(M)}=\Bgnorm{\sum_{k=1}^n
e_{k1}\ot x_k}_{S^p(L^p(M))}.
\end{equation}
The space $L^p(M,\ell^2_c)$ is the completion of $L^p(M)\ot \ell^2$
for this norm. It identifies to the space of sequences $(x_k)_{k\geq
1}$ in $L^p(M)$ such that $\sum_{k=1}^{+\infty} x_k\ot e_k $ is
convergent for the above norm. We define $L^p(M,\ell^2_r)$
similarly. For any finite family $(x_k)_{1 \leq k\leq n}$ in
$L^{p}(M)$, we have
\begin{equation*}
\Bgnorm{\sum_{k=1}^n x_k\ot
e_k}_{L^p(M,\ell^2_r)}=\Bgnorm{\bigg(\sum_{k=1}^n|x_k^*|^2\bigg)^{\frac{1}{2}}}_{L^p(M)}=\Bgnorm{\sum_{k=1}^ne_{1k}\ot
x_k}_{S^p(L^p(M))}.
\end{equation*}
For any $1\leq  p <\infty$ and for any $x_1,\ldots,x_n\in L^p(M)$,
we have
\begin{equation}\label{formule dualite}
\Bgnorm{\sum_{k=1}^n x_k\ot
e_k}_{L^p(M,\ell^2_c)}=\sup\Bigg\{\left|\sum_{k=1}^n\langle
x_k,y_k\rangle_{L^p(M),L^{p^*}(M)}\right|\ \colon \
\Bgnorm{\sum_{k=1}^n y_k\ot e_k}_{L^{p^*}(M,\ell^2_r)}\leq 1\Bigg\}.
\end{equation}
A similar formula holds for the space $L^p(M,\ell^2_r)$. For
simplicity, we write $S^p(\ell^2_c)$ for
$L^p\big(B(\ell^2),\ell^2_c\big)$. If $2 \leq p<\infty$ we define
the Banach space $L^p(M,\ell^2_\rad)=L^p(M,\ell^2_c)\cap
L^p(M,\ell^2_r)$. For any $u\in L^p(M,\ell^2_\rad)$, we have
\begin{equation*}
\norm{u}_{L^p(\ell^2_\rad)}=\max\Big\{
\norm{u}_{L^p(M,\ell^2_c)},\norm{u}_{L^p(M,\ell^2_r)}\Big\}.
\end{equation*}
If $1 \leq p \leq 2$ we define the Banach space
$L^p(M,\ell^2_\rad)=L^p(M,\ell^2_c)+ L^p(M,\ell^2_r)$. For any $u\in
L^p(M,\ell^2_\rad)$, we have
\begin{equation*}
\norm{u}_{L^p(M,\ell^2_\rad)}=\inf\Big\{
\norm{u_1}_{L^p(M,\ell^2_c)}+\norm{u_2}_{L^p(M,\ell^2_r)}\Big\}.
\end{equation*}
where the infimum runs over all possible decompositions $u=u_1+u_2$
with $u_1 \in L^p(M,\ell^2_c)$ and $u_2\in L^p(M,\ell^2_r)$. Recall
that, if $1<p<\infty$, we have an isometric identification
\begin{equation}
\label{dualite Lp(M,ell2rad)}
 L^p(M,\ell^2_\rad)^*=L^{p^*}(M,\ell^2_\rad).
\end{equation}
Let $X$ be a Banach space and let $(\varepsilon_k)_{k\geq 1}$ be a
sequence of independent Rademacher variables on some probability
space $\Omega$. Let $\Rad(X)\subset L^2(\Omega;X)$ be the closure of
${\rm Span}\big\{\varepsilon_k\otimes x\, :\, k\geq 1,\ x\in
X\big\}$ in the Bochner space $L^2(\Omega;X)$. Thus for any finite
family $x_1,\ldots,x_n$ in $X$, we have
$$
\Bgnorm{\sum_{k=1}^{n} \varepsilon_k\otimes x_k}_{\Rad(X)} \,=\,
\Bigg(\int_{\Omega}\bgnorm{\sum_{k=1}^{n} \varepsilon_k(\omega)\,
x_k}_{X}^{2}\,d\omega\,\Bigg)^{\frac{1}{2}}.
$$
If $1\leq p<\infty$, the noncommutative Khintchine's inequalities
(see \cite{LPP} and \cite{PX}) implies
\begin{equation}
\label{noncommutative Khintchine} \Rad\big(L^p(M)\big)\approx
L^p(M,\ell^2_\rad).
\end{equation}
We say that a set $\mathcal{F} \subset B(X)$ is $R$-bounded if there
is a constant $C\geq 0$ such that for any finite families
$T_1,\ldots, T_n$ in $\mathcal{F}$, and $x_1,\ldots,x_n$ in $X$, we
have
$$
\Bgnorm{\sum_{k=1}^{n} \varepsilon_k\otimes T_k
(x_k)}_{\Rad(X)}\,\leq\, C\, \Bgnorm{\sum_{k=1}^{n}
\varepsilon_k\otimes x_k}_{\Rad(X)}.
$$
In this case, we let $R(\mathcal{F})$ denote the smallest possible
$C$, which is called the $R$-bound of $\mathcal{F}$. $R$-boundedness
was introduced in \cite{BG} and then developed in the fundamental
paper \cite{ClP}. We refer to the latter paper and to \cite[Section
2]{KW} for a detailed presentation.

On noncommutative $L^p$-spaces, it will be convenient to consider
two naturals variants of this notion, introduced in \cite[Chapter
4]{JMX}. Let $1 < p<\infty$. A subset $\mathcal{F}$ of
$B\big(L^p(M)\big)$ is Col-bounded (resp. Row-bounded) if there
exists a constant $C\geq 0$ such that for any finite families
$T_1,\ldots,T_n$ in $\mathcal{F}$, and $x_1,\ldots,x_n$ in $L^p(M)$,
we have
\begin{equation}
\label{colbound}
\Bgnorm{\bigg(\sum_{k=1}^{n}\big|T_k(x_k)\big|^2\bigg)^{\frac{1}{2}}}_{L^p(M)}\leq
C\Bgnorm{\bigg(\sum_{k=1}^{n}|x_k|^2\bigg)^{\frac{1}{2}}}_{L^p(M)}
\end{equation}

\begin{equation}
\label{rowbound} \Bigg( \text{resp.
}\Bgnorm{\bigg(\sum_{k=1}^{n}\big|T_k(x_k)^*\big|^2\bigg)^{\frac{1}{2}}}_{L^p(M)}\leq
C\Bgnorm{\bigg(\sum_{k=1}^{n}|x_k^*|^2\bigg)^{\frac{1}{2}}}_{L^p(M)}\Bigg).
\end{equation}
The least constant $C$ satisfying (\ref{colbound}) will be denoted 
by $\col(\mathcal{F})$. Obviously any Rad-bounded (resp.
Col-bounded, resp. Row-bounded) set is bounded. It follows from
(\ref{noncommutative Khintchine}) that if a subset $\mathcal{F}$ of
$B\big(L^p(M)\big)$ is both Col-bounded and Row-bounded, then it is
Rad-bounded.

Note that contrary to the case of $R$-boundedness, a singleton
$\{T\}$ is not automatically Col-bounded or Row-bounded. Indeed,
$\{T\}$ is Col-bounded (resp. Row-bounded)  if and only if $T \ot
I_{\ell^2}$ extends to a bounded operator on $L^p(M,\ell^2_c)$
(resp. $L^p(M,\ell^2_r)$). And it turns out that if
$1<p\not=2<\infty$, according to \cite[Example 4.1]{JMX}, there
exists a bounded operator $T$ on $S^p$ such that $T\ot I_{\ell^2}$
does not extend to a bounded operator on $S^p(\ell^2_c)$. Moreover,
$T\ot I_{\ell^2}$ extends to a bounded operator on $S^p(\ell^2_r)$.
Then, we also deduce that there are sets $\mathcal{F}$ which are
Rad-bounded and Col-bounded without being Row-bounded. Similarly,
one may find sets which are Rad-bounded and Row-bounded without
being Col-bounded, or which are Rad-bounded without being either
Row-bounded or Col-bounded.

%

We turn to Ritt operators, the key class of this paper, and recall
some of their main features. Details and complements can be found in
\cite{ALM}, \cite{Bl1}, \cite{Bl2}, \cite{LM0}, \cite{Ly},
\cite{NZ}, \cite{N} and \cite{V}. Let $X$ be a Banach space. We say
that an operator $T\in B(X)$ is a Ritt operator if the two sets
\begin{equation}
\label{ensembles de Ritt} \big\{T^n\,:\, n\geq
0\big\}\qquad\hbox{and}\qquad \big\{n(T^n-T^{n-1}) \,:\, n\geq
1\big\}
\end{equation}
are bounded. This is equivalent to the spectral inclusion
\begin{equation}\label{spectral inclusion}
\sigma(T)\subset \ovl{\mathbb{D}}
\end{equation}
and the boundedness of the set
\begin{equation}\label{ensemble de Ritt}
\bigl\{(\lambda-1)R(\lambda,T)\,:\,\vert\lambda\vert>1\bigr\}
\end{equation}
where $R(\lambda,T)=(\lambda I-T)^{-1}$ denotes the resolvent
operator and $\mathbb{D}$ denotes the open unit disc centered at 0.
Likewise we say that $T$ is an $R$-Ritt operator if the two sets in
(\ref{ensembles de Ritt}) are $R$-bounded. This is equivalent to the
inclusion (\ref{spectral inclusion}) and the $R$-boundedness of the
set (\ref{ensemble de Ritt}).

Let $T$ be a Ritt operator. The boundedness of (\ref{ensemble de
Ritt}) implies the existence of a constant $K\geq 0$ such that
$|\lambda-1|\bnorm{R(\lambda,T)}_{X \to X}\leq K$ whenever
$\Re(\lambda)>1$. This means that $I-T$ is a sectorial operator.
Thus for any $\alpha > 0$, one can consider the fractional power
$(I-T)^\alpha$. We refer to \cite[Chapter 3]{Haa}, \cite{KW} and
\cite{MCS} for various definitions of these (bounded) operators and
their basic properties.

We will use the following two naturals variants of the notion of
$R$-Ritt operator.
\begin{defi}
Suppose $1<p<\infty$. Let $T$ be a bounded operator on $L^p(M)$. We
say that $T$ is a Col-Ritt (resp. Row-Ritt) operator if the two sets
(\ref{ensembles de Ritt}) are Col-bounded (resp. Row-bounded).
\end{defi}

\begin{remark}
Assume that $1<p<\infty$. Let $T$ be a bounded operator on $L^p(M)$.
Using (\ref{formule dualite}), it is easy to see that $T$ is
Col-Ritt if and only if $T^*$ is Row-Ritt on $L^{p^*}(M)$.
\end{remark}

We let $\mathcal{P}$ denote the algebra of all complex polynomials.
Let $T$ be a bounded operator on a Banach space $X$. Let $\gamma \in
\big]0,\frac{\pi}{2}\big[$. Accordingly with \cite{LM0}, we say that
$T$ has a bounded $H^\infty(B_{\gamma})$ functional calculus if and
only if there exists a constant $K\geq 1$ such that
$$
\bnorm{\varphi(T)}_{X\to X}\leq
K\norm{\varphi}_{H^\infty(B_{\gamma})}
$$
for any $\varphi\in \mathcal{P}$. Naturally, we let:

\begin{defi}
Suppose $1< p<\infty$. Let $T$ be a bounded operator on $L^p(M)$.
Let $\gamma \in \big]0,\frac{\pi}{2}\big[$. We say that $T$ admits a
completely bounded $H^\infty(B_{\gamma})$ functional calculus if $T$
is completely bounded and if $I_{S^p}\ot T$ admits a bounded
$H^\infty(B_{\gamma})$ functional calculus on $S^p\big(L^p(M)\big)$.
\end{defi}
Let $T$ be a bounded operator on $L^p(M)$ and $\gamma \in
\big]0,\frac{\pi}{2}\big[$. Note that $T$ admits a completely
bounded $H^\infty(B_{\gamma})$ functional calculus if and only if
there exists a constant $K\geq 1$ such that
$$
\bnorm{\varphi(T)}_{cb, L^p(M)\to L^p(M)}\leq
K\norm{\varphi}_{H^\infty(B_{\gamma})}
$$
for any $\varphi\in \mathcal{P}$.

\section{Results related to Col-Ritt or Row-Ritt operators}


%
%
%
In the subsequent sections, we need some preliminary results on
Col-Ritt or Row-Ritt operators that we present here. Some of them
are analogues of existing results in the context of $R$-Ritt
operators, for which we will omit proofs.

We start with a variant of \cite[Proposition 2.8]{ALM} suitable with
our context. The proof is similar, using \cite[Lemma 4.2]{JMX}
instead of \cite[Lemma 2.1]{ALM}.

\begin{prop}
\label{prop bornitude} Suppose $1<p<\infty$. Let $T$ be a Col-Ritt
operator on $L^p(M)$. For any $\alpha>0$, the set
$$
\Big\{n^\alpha (\varrho T)^{n-1}(I-\varrho T)^{\alpha}\, \colon\,
n\geq 1,\ \varrho\in]0,1]\Big\}
$$
is Col-bounded. Moreover, a similar result holds for Row-Ritt
operators.
\end{prop}


Moreover, we need the following result \cite{LM0}.
\begin{thm}
\label{LM FC implies R-Ritt} Suppose $1<p<\infty$. Let $T$ be a
bounded operator on $L^p(M)$ with a bounded $H^\infty(B_\gamma)$
functional calculus for some $\gamma\in \big]0,\frac{\pi}{2}\big[$.
Then $T$ is R-Ritt.
\end{thm}

In the next statement, we establish a variant of the above result.

\begin{thm}
\label{FC CB implies Col} Suppose $1<p<\infty$. Let $T$ be a bounded
operator on $L^p(M)$. Assume that $T$ admits a completely bounded
$H^\infty(B_\gamma)$ functional calculus for some $\gamma\in
\big]0,\frac{\pi}{2}\big[$. Then the operator $T$ is both Col-Ritt
and Row-Ritt.
\end{thm}

\begin{preuve}
We will only show the `column' result, the proof for the `row' one
being the same. We wish to show that the sets
$$
\mathcal{F}=\big\{T^m\,:\, m\geq 0\big\}\qquad\hbox{and}\qquad
\mathcal{G}=\big\{m(T^m-T^{m-1}) \,:\, m\geq 1\big\}
$$
are Col-bounded. We consider the operator $I\ot T$ on the
noncommutative $L^p$-space $S^p\big(L^p(M)\big)$. Then, applying
Theorem \ref{LM FC implies R-Ritt}, we obtain that the sets
$$
\mathcal{T}=\big\{I_{S^p}\ot T^m\,:\, m\geq
0\big\}\qquad\hbox{and}\qquad \mathcal{K}=\big\{mI_{S^p}\ot
(T^m-T^{m-1}) \,:\, m\geq 1\big\}
$$
are Rad-bounded. Now consider $x_1,\ldots,x_n$ in $L^p(M)$ and
$T_1,\ldots,T_n$ in $\mathcal{F}$. For any finite sequence
$(\varepsilon_k)_{1\leq k\leq n}$ valued in $\{-1,1\}$, we have
\begin{align*}
\Bgnorm{\bigg(\sum_{k=1}^{n}|x_k|^2\bigg)^{\frac{1}{2}}}_{L^p(M)}
    &=\Bgnorm{\bigg(\sum_{k=1}^{n}(\varepsilon_k x_k)^*(\varepsilon_k x_k)\bigg)^{\frac{1}{2}}}_{L^p(M)}\\
    &=\Bgnorm{\sum_{k=1}^{n}\varepsilon_k e_{k1}\ot x_k}_{S^p(L^p(M))}.
\end{align*}
Then passing to the average over all possible choices of
$\varepsilon_k=\pm 1$, we obtain that
$$
\Bgnorm{\bigg(\sum_{k=1}^{n}|x_k|^2\bigg)^{\frac{1}{2}}}_{L^p(M)}=
\Bgnorm{\sum_{k=1}^{n}\varepsilon_k\ot e_{k1}\ot
x_k}_{\Rad(S^p(L^p(M)))}.
$$
By a similar computation, we have
$$
\Bgnorm{\bigg(\sum_{k=1}^{n}|T_k(x_k)|^2\bigg)^{\frac{1}{2}}}_{L^p(M)}
=\Bgnorm{\sum_{k=1}^{n}\varepsilon_k\ot (I_{S^p} \ot T_k)(e_{k1}\ot
x_k)}_{\Rad(S^p(L^p(M)))}.
$$
It follows that
$$
\Bgnorm{\bigg(\sum_{k=1}^{n}|T_k(x_k)|^2\bigg)^{\frac{1}{2}}}_{L^p(M)}\leq
\Rad(\mathcal{T})\Bgnorm{\bigg(\sum_{k=1}^{n}|x_k|^2\bigg)^{\frac{1}{2}}}_{L^p(M)}.
$$
This concludes the proof of Col-boundedness of $\mathcal{F}$ with
$\col(\mathcal{F}) \leq \Rad(\mathcal{T})$. The proof for the set
$\mathcal{G}$ is identical.
\end{preuve}

\begin{remark}
Suppose $1<p\not=2<\infty$. The complete boundedness assumption in
Theorem \ref{FC CB implies Col} cannot be replaced by a boundedness
assumption.
\end{remark}

\begin{preuve}
We have already recalled that, there exists a bounded operator $T$
on $S^p$ such that $\{T\}$ is not Col-bounded. Let us fix $\gamma
\in \big]0,\frac{\pi}{2}\big[$. We may clearly assume that
$\sigma(T)$ is included in the open set $B_\gamma$. Using the
Dunford calculus, it is easy to prove that $T$ is a Ritt operator
which admits a bounded $H^\infty(B_\gamma)$ functional calculus. The
set $\{T\}$ is not Col-bounded. Hence $T$ cannot be Col-Ritt.
\end{preuve}

Now, we give a precise definition of `square functions' which
clarifies (\ref{norm T1}), (\ref{norm T1 prime}) and (\ref{column
and row square functions}) and a few comments. Let $T$ a Ritt
operator on $L^p(M)$. For any $\alpha>0$, let us consider
$$
x_k=k^{\alpha-\frac{1}{2}}T^{k-1}(I-T)^{\alpha}(x)
$$
for any $k\geq 1$. If the sequence belongs to the space
$L^p(M,\ell^2_{c})$, then $\norm{x}_{p,T,c,\alpha}$ is defined as
the norm of $(x_k)_{k\geq 1}$ in that space. Otherwise, we set
$\norm{x}_{p,T,c,\alpha}=\infty$. In particular,
$\norm{x}_{p,T,c,\alpha}$ can be infinite. We define the quantities
$\norm{x}_{p,T,r,\alpha}$ by the same way. The quantities
$\norm{x}_{p,T,\alpha}$ are defined similarly in \cite{ALM}, using
the space $L^p(M,\ell^2_{\rad})$ instead of $L^p(M,\ell^2_{c})$.

%
%

Finally, note that, if $2\leq p < \infty$, we have
$$
\norm{x}_{p,T,\alpha}=\max\big\{\norm{x}_{p,T,c,\alpha},\norm{x}_{p,T,r,\alpha}\big\}.
$$
and if $1\leq p\leq 2$, we have
$$
\norm{x}_{p,T,\alpha}=\inf\Big\{\norm{u}_{L^{p}(M,\ell^2_c)}+\norm{v}_{L^{p}(M,\ell^2_r)}\
\colon\  u_k+v_k=k^{\alpha-\frac{1}{2}}T^{k-1}(I-T)^{\alpha}x\text{
for any integer $k$}\Big\}.
$$
In \cite{LM0}, the following connection between the boundedness of
square functions and functional calculus is established.
\begin{thm}
\label{principal thm LM ghost} Suppose $1<p<\infty$. Let $T$ be a
bounded operator on $L^p(M)$. The following assertions are
equivalent.

\begin{enumerate}
  \item The operator $T$ is R-Ritt and $T$ and its adjoint $T^*$ both satisfy uniform
  estimates
  $$
\norm{x}_{p,T,1} \lesssim \norm{x}_{L^{p}(M)} \ \ \ \ \ \text{and} \
\ \ \ \ \norm{y}_{p^*,T^*,1} \lesssim\norm{y}_{L^{p^*}(M)}
  $$
  for any $x\in L^p(M)$ and $y\in L^{p^*}(M)$.
  \item  The operator $T$ admits a bounded $H^\infty(B_\gamma)$ functional
calculus for some $\gamma\in \big]0,\frac{\pi}{2}\big[$.
\end{enumerate}
\end{thm}
Recall a special case of the principal result of \cite{ALM}.
\begin{thm}
\label{thm principal ALM} Let $T$ be an $R$-Ritt operator on
$L^p(M)$ with $1<p<\infty$. For any $\alpha,\beta>0$ we
  have an equivalence
  $$
  \norm{x}_{p,T,\alpha} \approx \norm{x}_{p,T,\beta},\qquad x\in L^p(M).
  $$
\end{thm}
We shall now present a variant suitable to our context.

For any integer $n\geq 1$, we identify the algebra $M_n$ of all
$n\times n$ matrices with the space of linear maps $\ell^2_n\to
\ell^2_n$. 
For any infinite matrix $[c_{ij}]_{ i,j\geq 1}$,
we set
$$
\bnorm{[c_{ij}]}_{\rm reg}=\sup_{n\geq 1}\Bnorm{\bigr[\vert
c_{ij}\vert\,\bigr]_{1\leq i,j\leq n}}_{B(\ell^2_n)}
$$
This is the so-called `regular norm'. We refer to \cite{Pis11} and
\cite{Pis5} for more information on regular norms.


The next proposition will be useful. This result is similar to
\cite[Proposition 2.6]{ALM}.

\begin{prop}
\label{prop 2.5 version col} Suppose $1<p<\infty$. Let
$[c_{ij}]_{i,j\geq 1}$ be an infinite matrix with
$\bnorm{[c_{ij}]}_{\rm reg}<\infty$. Suppose that $\big\{T_{ij}\
\colon\ i,j \geq 1 \big\}$ is a Col-bounded set of operators on
$L^p(M)$. Then the linear map
$$
\begin{array}{cccc}
 \big[c_{ij}T_{ij}\big]:   &  L^p\big(M,\ell^2_c\big)   &  \longrightarrow   &  L^p\big(M,\ell^2_c\big)  \\
    & \displaystyle \sum_{j=1}^{+\infty} x_j \ot e_j  &  \longmapsto       & \displaystyle \sum_{i=1}^{+\infty} \Bigg(\sum_{j=1}^{+\infty}c_{ij}T_{ij}(x_j)\Bigg)\ot e_i  \\
\end{array}
$$
is well-defined and bounded. Moreover, we have a similar result for
Row-bounded sets.
\end{prop}

\begin{preuve}
We shall only prove the `Col' result. We can assume that
$\bnorm{[c_{ij}]}_{\rm reg} \leq 1$. Let $n\geq 1$. By \cite[Lemma
2.2]{ALM}, we can write $c_{ij}=a_{ij}b_{ij}$ for any $1 \leq i,j
\leq n$ with
$$
\sup_{1 \leq i\leq n}\sum_{j=1}^{n}|a_{ij}|^2\leq 1 \ \ \ \ \
\text{and} \ \ \ \ \ \sup_{1 \leq j \leq
n}\sum_{i=1}^{n}|b_{ij}|^2\leq 1.
$$
Let $x_1,\ldots,x_n\in L^p(M)$ and $y_1,\ldots,y_n \in L^{p^*}(M)$.
Since the set $\big\{T_{ij}\ |\ i,j \geq 1 \big\}$ is Col-bounded,
there exists a positive constant $C$ such that
\begin{align*}
\MoveEqLeft
 \left|\sum_{i=1}^n \Bigg\langle \sum_{j=1}^nc_{ij}T_{ij} (x_j),y_i\Bigg\rangle_{L^p(M),L^{p^*}(M)}\right|
 =  \left|\sum_{i,j=1}^{n} \big\langle a_{ij}b_{ij} T_{ij}(x_j),y_i\big\rangle_{L^p(M),L^{p^*}(M)}\right|\\
 &=  \left|\sum_{i,j=1}^{n} \big\langle T_{ij}(b_{ij}x_j),a_{ij}y_i\big\rangle_{L^p(M),L^{p^*}(M)}\right|\\
 &\leq
 \Bgnorm{\bigg(\sum_{i,j=1}^{n} |T_{ij}(b_{ij}x_j)|^2\bigg)^{\frac{1}{2}}}_{L^{p}(M)}
 \Bgnorm{\bigg(\sum_{i,j=1}^{n}\big|(a_{ij}y_i)^*\big|^2\bigg)^{\frac{1}{2}}}_{L^{p^*}(M)}\\
 &\leq C
 \Bgnorm{\bigg(\sum_{i,j=1}^{n}|b_{ij}x_j|^2\bigg)^{\frac{1}{2}}}_{L^{p}(M)}
 \Bgnorm{\bigg(\sum_{i,j=1}^{n}|a_{ij}y_i^*|^2\bigg)^{\frac{1}{2}}}_{L^{p^*}(M)}.
\end{align*}
Now, we have
$$
 \sum_{i,j=1}^{n}|b_{ij}x_j|^2=  \sum_{j=1}^{n}|x_j|^2\Bigg(\sum_{i=1}^{n}|b_{ij}|^2\Bigg)\leq  \sum_{j=1}^{n}|x_j|^2.
$$
Similarly, we have
\begin{align*}
 \sum_{i,j=1}^{n}| a_{ij}y_i^*|^2
   &\leq   \sum_{i=1}^{n}|y_i^*|^2.
\end{align*}
Consequently
$$
\left|\sum_{i=1}^n \Bigg\langle\sum_{j=1}^n c_{ij}T_{ij}(x_j),y_i
\Bigg\rangle_{L^p(M),L^{p^*}(M)}\right| \leq
C\Bgnorm{\bigg(\sum_{j=1}^{n}|x_j|^2
\bigg)^{\frac{1}{2}}}_{L^{p}(M)}\Bgnorm{\bigg(
\sum_{i=1}^{n}|y_i^*|^2\bigg)^{\frac{1}{2}}}_{L^{p^*}(M)}.
$$
Taking the supremum over all $y_1,\ldots,y_n \in L^{p^*}(M)$ such
that $\bnorm{(\sum_{i=1}^{n}|y_i^*|^2
)^{\frac{1}{2}}}_{L^{p^*}(M)}\leq 1$, we obtain
$$
\Bgnorm{\sum_{i=1}^n \bigg(\sum_{j=1}^n c_{ij}T_{ij}(x_j)\bigg) \ot
e_i}_{L^{p}(M,\ell^2_c)} \leq C\Bgnorm{\sum_{j=1}^{n}x_j \ot e_j
}_{L^{p}(M,\ell^2_c)}
$$
by (\ref{formule dualite}). We conclude with \cite[Corollary
2.12]{JMX}.
\end{preuve}

Now, we state a result which allows to estimate square functions
$\norm{\cdot}_{p,T,c,\alpha}$ and $\norm{\cdot}_{p,T,r,\alpha}$ by
means
of approximation processes, whose proof is similar to \cite[Lemma 3.2]{ALM}.
\begin{lemma}
\label{lemma approx} Suppose $1<p<\infty$. Assume that $T$ is a
Col-Ritt operator on $L^p(M)$. Let $\alpha>0$.
\begin{enumerate}
  \item Let $V$ be an operator on $L^p(M)$ such that $TV=VT$ with $\{V\}$ Col-bounded.
  Then, for any $x\in L^p(M)$, we have
$$
\norm{V(x)}_{p,T,c,\alpha}\leq
\col\big(\{V\}\big)\norm{x}_{p,T,c,\alpha}.
$$
  \item Let $\nu\geq \alpha +1$ be an integer and let $x\in
\Ran\big((I-T)^\nu\big)$. Then
$$
\norm{x}_{p,\varrho T,c,\alpha}\xra[\varrho \to
1^-]{}\norm{x}_{p,T,c,\alpha}.
$$
\end{enumerate}
Moreover, the same result holds with $\norm{\cdot}_{p,T,c,\alpha}$
replaced by $\norm{\cdot}_{p,T,r,\alpha}$ for Row-Ritt operators.
\end{lemma}

Now we state an equivalence result in our context similar to Theorem
\ref{thm principal ALM}.
\begin{thm}
\label{Thm equivalence square function colrow} Let $T$ be a bounded
operator on $L^p(M)$ with $1<p<\infty$. Let $\alpha,\beta>0$.
\begin{enumerate}
  \item If $T$ is Col-Ritt, we
  have an equivalence
  $$
  \norm{x}_{p,T,c,\alpha} \approx \norm{x}_{p,T,c,\beta},\ \ \ \ \ \ x\in L^p(M).
  $$
  \item If $T$ is Row-Ritt, we
  have an equivalence
  $$
  \norm{x}_{p,T,r,\alpha} \approx \norm{x}_{p,T,r,\beta},\ \ \ \ \ \ x\in L^p(M).
  $$
\end{enumerate}
\end{thm}

\begin{preuve}
The proof is similar to the one of \cite[Theorem 3.3]{ALM}, using
Proposition \ref{prop bornitude}, Proposition \ref{prop 2.5 version
col}, Lemma \ref{lemma approx} and \cite[Corollary 2.12]{JMX}.
\end{preuve}

\section{Comparison between squares functions and the usual norm}

We aim at showing Theorem \ref{Th 3}. We will provide an example on
the Schatten space $S^p$. This example also prove that in general,
row and column square functions are not equivalent (Theorem \ref{thm
non equivalence}).

Let $a$ a bounded operator on $\ell^2$. Assume $1<p<\infty$. We let
$\mathcal{L}_a \colon S^p\to S^p$ the left multiplication by $a$ on
$S^p$ defined by $\mathcal{L}_a(x)=ax$ and we denote $\mathcal{R}_a
\colon S^p\to S^p$ the right multiplication. It is clear that
$\mathcal{L}_a^*$ and $\mathcal{R}_a^*$ are the right multiplication
and the left multiplication by $a$ on $S^{p^*}$. Note that, by
\cite[Proposition 8.4 (4)]{JMX}, if $I-a$ has dense range then
$\Ran(I-\mathcal{L}_a)$ is dense in $S^p$. The next statement gives
a link between properties of $a$ and its associated multiplication
operators.

\begin{prop}
\label{Ritt and La} Suppose $1<p<\infty$. Assume that $a$ is a
bounded operator on $\ell^2$.
\begin{enumerate}
  \item  If $a$ is a Ritt operator then the left multiplication $\mathcal{L}_a$ is a Ritt operator on $S^p$.
  \item  Let $\gamma\in \big]0,\frac{\pi}{2}\big[$. Then $\mathcal{L}_a$ has a bounded $H^\infty(B_\gamma)$
         functional calculus if and only if $a$ has one. In that case, $\mathcal{L}_a$
         actually has a completely bounded $H^\infty(B_\gamma)$ functional calculus.
\end{enumerate}
Moreover, we have a similar result for right multiplication.
\end{prop}
%

\begin{preuve}
We have $\sigma(\mathcal{L}_a)\subset \sigma(a)$. Moreover, if
$\lambda\in \rho(a)$ we have
$R(\lambda,\mathcal{L}_a)=\mathcal{L}_{R(\lambda,a)}$. The first
assertion clearly follows. The statement (2) is a straightforward
consequence of
$$
I_{S^p} \ot \mathcal{L}_a=\mathcal{L}_{I_{\ell^2}\ot a} \ \ \ \
\text{and}\ \ \ \ f(\mathcal{L}_a) = \mathcal{L}_{f(a)},\qquad f\in
\mathcal{P}.
$$
The proof of the `right' result is identical.
\end{preuve}

%
We denote by $(e_k)_{k\geq 1}$ the canonical basis of $\ell^2$. Now,
for any integer $k\geq 1$, we fix $a_k=1-\frac{1}{2^k}$. We consider
the selfadjoint bounded diagonal operator $a$ on $\ell^2$ defined by
\begin{equation}\label{def de a}
a\Bigg(\sum_{k=1}^{+\infty}x_ke_k\Bigg)=\sum_{k=1}^{+\infty}a_kx_ke_k.
\end{equation}
It follows from the Spectral Theorem for normal operators, that the
operator $a$ admits a bounded $H^\infty(B_\gamma)$ functional
calculus for any $\gamma\in \big]0,\frac{\pi}{2}\big[$. Thus
$\mathcal{L}_a$ and $\mathcal{R}_a$ admit a completely bounded
$H^\infty(B_\gamma)$ functional calculus for any $\gamma\in
\big]0,\frac{\pi}{2}\big[$ (hence $\mathcal{L}_a$ and
$\mathcal{R}_a$ are Ritt operators).
\begin{lemma}
\label{equivalence usual norm and col for La} Assume that $2 \leq
p<\infty$. Let $a$ be the bounded operator on $\ell^2$ defined by
(\ref{def de a}). If $\mathcal{L}_a \colon S^p\to S^p$ and
$\mathcal{R}_a \colon S^p\to S^p$ are the left and right
multiplication operators associated to $a$, we have
\begin{equation}
\label{equivalence column square function}
    \norm{x}_{p,\mathcal{L}_a,c,1} \approx \norm{x}_{S^p}
    \ \ \ \ \ \text{and} \ \ \ \ \    \norm{x}_{p,\mathcal{R}_a,r,1} \approx
    \norm{x}_{S^p}, \qquad x\in S^p.
\end{equation}
\end{lemma}

\begin{preuve}
We will only show the result for the operator $\mathcal{L}_a$, the
proof for $\mathcal{R}_a$ being the same. For any $x\in S^p$ and any
$\varrho \in ]0,1[$, we have
\begin{align*}
k\big((\varrho\mathcal{L}_{a})^{k-1}(I-\varrho\mathcal{L}_{a})(x)\big)^{*}\big((\varrho\mathcal{L}_{a})^{k-1}(I-\varrho\mathcal{L}_{a})(x)\big)
    &=k\big((\varrho a)^{k-1}(I-\varrho a)x\big)^{*}\big((\varrho a)^{k-1}(I-\varrho a)x\big)\\
    &= kx^*(I-\varrho a)(\varrho a)^{2(k-1)}(I-\varrho a)x\\
    &= kx^*(I-\varrho \mathcal{L}_{a})^2(\varrho\mathcal{L}_{a})^{2(k-1)}(x).
\end{align*}
Now, for any $z\in \mathbb{D}$, we have
\begin{equation}
\label{power series}
    \sum_{k=1}^{+\infty} kz^{k-1}=(1-z)^{-2}.
\end{equation}
Since the operator $\mathcal{L}_{a}$ is a contraction, we deduce
that, for every $\varrho \in ]0,1[$, the operator $I-(\varrho
\mathcal{L}_{a})^2$ is invertible and that we have
\begin{equation}\label{(I-a)^-2=}
   \sum_{k=1}^{+\infty}k (\varrho \mathcal{L}_{a})^{2(k-1)}=\big(I-(\varrho
   \mathcal{L}_{a})^2\big)^{-2},
\end{equation}
the series being absolutely convergent. Then we deduce that the
series
$$
\sum_{k=1}^{+\infty}k\big((\varrho\mathcal{L}_{a})^{k-1}(I-\varrho\mathcal{L}_{a})(x)\big)^{*}\big((\varrho\mathcal{L}_{a})^{k-1}(I-\varrho\mathcal{L}_{a})(x)\big)
$$
is convergent in the Banach space $S^{\frac{p}{2}}$ and that
\begin{align*}
\sum_{k=1}^{+\infty}k\big((\varrho\mathcal{L}_{a})^{k-1}(I-\varrho\mathcal{L}_{a})(x)\big)^{*}\big((\varrho\mathcal{L}_{a})^{k-1}(I-\varrho\mathcal{L}_{a})(x)\big)
    &=x^*(I-\varrho \mathcal{L}_{a})^2\big(I-(\varrho \mathcal{L}_{a})^2\big)^{-2}x\\
    &=x^*(I+\varrho a)^{-2}x.
\end{align*}
We deduce that
\begin{align*}
\norm{x}_{\varrho\mathcal{L}_a,c,1}
&=\bgnorm{\Big(x^*(I+\varrho a)^{-2}x\Big)^{\frac{1}{2}}}_{S^p}\\
&=\bnorm{(I+\varrho a)^{-1}x}_{S^p}.
\end{align*}
Then, for any $x\in S^p$, we obtain the estimate
\begin{align*}
\norm{x}_{p,\varrho\mathcal{L}_a,c,1} &\leq \bnorm{(I+\varrho
a)^{-1}}_{B(\ell^2)}\norm{x}_{S^p}\\
&\leq\norm{x}_{S^p}.
\end{align*}
By a similar computation, for any $x\in S^p$, we have
$$
\frac{1}{2}\norm{x}_{S^p}\leq\norm{x}_{p,\varrho\mathcal{L}_a,c,1}.
$$
%
Applying Lemma \ref{lemma approx} (2), we deduce an equivalence
$$
\frac{1}{2}\norm{x}_{S^p}\leq\norm{x}_{p,\mathcal{L}_a,c,1}\leq
\norm{x}_{S^p},\qquad x \in \Ran\big((I-\mathcal{L}_a)^2\big).
$$
For any integer $n\geq 1$, we let $d_n$ the bounded diagonal
operator on $\ell^2$ defined by the matrix
$\diag(1,\ldots,1,0,\ldots)$. It is not difficult to see that, for
any integer $n\geq 1$, the range of $\mathcal{L}_{d_n}$ is a
subspace of $\Ran\big((I-\mathcal{L}_a)^2\big)$. Hence we actually
have
\begin{equation*}
\frac{1}{2}\bnorm{\mathcal{L}_{d_n}(x)}_{S^p}\leq\bnorm{\mathcal{L}_{d_n}(x)}_{p,\mathcal{L}_a,c,1}\leq
\bnorm{\mathcal{L}_{d_n}(x)}_{S^p},\qquad x \in S^p,\quad n\geq 1.
\end{equation*}
Then, on the one hand, we obtain
\begin{equation*}
\bnorm{\mathcal{L}_{d_n}(x)}_{p,\mathcal{L}_a,c,1} \leq
\norm{x}_{S^p},\qquad x \in S^p,\quad n\geq 1.
\end{equation*}
By \cite[Corollary 2.12]{JMX} and (\ref{expression norm}), this
latter inequality is equivalent to
$$
\Bgnorm{\sum_{k=1}^{l}e_{k1} \ot
k^{\frac{1}{2}}\mathcal{L}_a^{k-1}(I-\mathcal{L}_a)\big(\mathcal{L}_{d_n}(x)\big)}_{S^p(S^p)}\lesssim
\norm{x}_{S^p},\qquad x \in S^p,\quad n\geq 1,\quad l\geq 1.
$$
Passing to the limit in the above inequality, we infer that
$$
\Bgnorm{\sum_{k=1}^{l}e_{k1} \ot
k^{\frac{1}{2}}\mathcal{L}_a^{k-1}(I-\mathcal{L}_a)(x)}_{S^p(S^p)}\lesssim
\norm{x}_{S^p},\qquad x \in S^p,\quad l\geq 1.
$$
Using again \cite[Corollary 2.12]{JMX}, we obtain that
\begin{equation*}
    \norm{x}_{p,\mathcal{L}_a,c,1} \leq \norm{x}_{S^p},\qquad x\in S^p.
\end{equation*}
Note, in particular that, for any $x\in S^p$, we have
$\norm{x}_{p,\mathcal{L}_a,c,1} <\infty$. On the other hand, note
that, for any integer $n\geq 1$, the operators $\mathcal{L}_a$ and
$\mathcal{L}_{d_n} $ commute. Hence, for any $x \in S^p$ and any
integer $n\geq 1$, we have
\begin{align*}
\norm{\mathcal{L}_{d_n}(x)}_{S^p}
  &\lesssim \bnorm{\mathcal{L}_{d_n}(x)}_{p,\mathcal{L}_a,c,1} \\
  &=\Bgnorm{\sum_{k=1}^{+\infty} e_{k1} \ot k^{\frac{1}{2}}\mathcal{L}_a^{k-1}(I-\mathcal{L}_a)\big(\mathcal{L}_{d_n}(x)\big) }_{S^p(S^p)}\\
  &=\Bgnorm{(I_{S^p}\ot \mathcal{L}_{d_n})\bigg(\sum_{k=1}^{+\infty} e_{k1} \ot k^{\frac{1}{2}} \mathcal{L}_a^{k-1}(I-\mathcal{L}_a)(x)\bigg)}_{S^p(S^p)}.
\end{align*}
Letting $n$ to the infinity, we deduce that
\begin{equation*}
\norm{x}_{S^p}  \lesssim  \norm{x}_{p,\mathcal{L}_a,c,1},\qquad x\in
S^p.
\end{equation*}
The proof is complete.
\end{preuve}

\begin{thm}
\label{thm non equivalence} Let $\alpha>0$. Let $a$ be the bounded
operator on $\ell^2$ defined by (\ref{def de a}). Let
$\mathcal{L}_a\colon S^p\to S^p$ and $\mathcal{R}_a\colon S^p\to
S^p$ be the left and right multiplication operators associated to
$a$. Assume that $2<p<\infty$. Then
\begin{equation}\label{c/r=infty}
\sup\Bigg\{\frac{\norm{x}_{p,\mathcal{L}_a,c,\alpha}}{\norm{x}_{p,\mathcal{L}_a,r,\alpha}}\
\colon\ x\in S^p\Bigg\}=\infty\ \ \ \ \ \text{and} \ \ \ \ \
\sup\Bigg\{\frac{\norm{x}_{p,\mathcal{R}_a,r,\alpha}}{\norm{x}_{p,\mathcal{R}_a,c,\alpha}}\
\colon\ x\in S^p\Bigg\}=\infty.
\end{equation}
Assume that $1<p<2$. Then
\begin{equation}\label{r/c=infty}
\sup\Bigg\{\frac{\norm{x}_{p,\mathcal{L}_a,r,\alpha}}{\norm{x}_{p,\mathcal{L}_a,c,\alpha}}\
\colon\ x\in S^p\Bigg\}=\infty \ \ \ \ \ \text{and} \ \ \ \ \
\sup\Bigg\{\frac{\norm{x}_{p,\mathcal{R}_a,c,\alpha}}{\norm{x}_{p,\mathcal{R}_a,r,\alpha}}\
\colon\ x\in S^p\Bigg\}=\infty.
\end{equation}
\end{thm}

\begin{preuve}
By Theorem \ref{Thm equivalence square function colrow}, it suffices
to prove the result for one specific real $\alpha$. Throughout the
proof, we will use $\alpha=1$. We first assume that $2<p<\infty$.
Given an integer $n\geq 1$, we consider $e=e_1+\cdots+e_n\in
\ell^2_n$ and $x=\frac{1}{\sqrt{n}}e\ot e\in S^p$. Clearly, we have
$$
xx^*=\sum_{i,j=1}^n e_{ij}.
$$
Now, we have
\begin{align*}
k\big(\mathcal{L}_{a}^{k-1}(I-\mathcal{L}_{a})(x)\big)\big(\mathcal{L}_{a}^{k-1}(I-\mathcal{L}_{a})(x)\big)^{*}
&=k\big(a^{k-1}(I-a)x\big)\big(a^{k-1}(I-a)x\big)^{*}\\
&=ka^{k-1}(I-a)xx^*(I-a)a^{k-1}\\
&=\sum_{i,j=1}^n ka^{k-1}(I-a)e_{ij}(I-a)a^{k-1}\\
&=\sum_{i,j=1}^n (1-a_i)(1-a_j)k(a_ia_j)^{k-1}e_{ij}.
\end{align*}
Using the equality (\ref{power series}), we obtain that the series
$$
\sum_{k=1}^{+\infty}
k\big(\mathcal{L}_{a}^{k-1}(I-\mathcal{L}_{a})(x)\big)\big(\mathcal{L}_{a}^{k-1}(I-\mathcal{L}_{a})(x)\big)^{*}
$$
is convergent in $S^{\frac{p}{2}}$ and that
$$
\sum_{k=1}^{+\infty}k\big(\mathcal{L}_{a}^{k-1}(I-\mathcal{L}_{a})(x)\big)\big(\mathcal{L}_{a}^{k-1}(I-\mathcal{L}_{a})(x)\big)^{*}
=\sum_{i,j=1}^n (1-a_i)(1-a_j)(1-a_ia_j)^{-2}e_{ij}.
$$
Now, note that
$$
(1-a_i)(1-a_j)(1-a_ia_j)^{-2}=\frac{2^{i+j}}{(2^i+2^j-1)^2}.
$$
We deduce that
\begin{align*}
\norm{x}_{p,\mathcal{L}_a,r,1}
&=\Bgnorm{\bigg(\sum_{i,j=1}^n\frac{2^{i+j}}{(2^i+2^j-1)^2}e_{ij}\bigg)^\frac{1}{2}}_{S^p}\\
&=\Bgnorm{\sum_{i,j=1}^n\frac{2^{i+j}}{(2^i+2^j-1)^2}e_{ij}}_{S^{\frac{p}{2}}}^\frac{1}{2}.
\end{align*}
We let $A=\Big[\frac{2^{i+j}}{(2^i+2^j-1)^2}\Big]_{1\leq i,j\leq n}$
be the $n\times n$ matrix in the last right member of the above
equations. We have
\begin{align*}
\norm{A}_{S^2_n}^2 &=\sum_{i,j=1}^{n}
\Bigg(\frac{2^{i+j}}{(2^i+2^j-1)^2}\Bigg)^2 \\
&=\sum_{i,j=1}^{n} \frac{4^{i+j}}{(2^i+2^j-1)^4}.
\end{align*}
Moreover, note that
\begin{align*}
 \frac{4^{i+j}}{(2^i+2^j-1)^4}
    &\leq 16 \frac{4^{i+j}}{(2^i+2^j)^4}\\
    &=16 \Bigg(\frac{1}{2^{i-j}+2^{j-i}+2}\Bigg)^2\\
    &\leq \frac{16}{4^{|i-j|}}.
\end{align*}
Thus we have
$$
\norm{A}_{S^2_n}^2\leq 32\bigg(\sum_{k\in
\mathbb{Z}}\frac{1}{4^{|k|}}\bigg)n\approx n.
$$

If $4\leq p<\infty$, we obtain
$$
\norm{x}_{p,\mathcal{L}_a,r,1}=\norm{A}_{S^{\frac{p}{2}}_n}^{\frac{1}{2}}\leq
\norm{A}_{S^2_n}^{\frac{1}{2}}\lesssim n^{\frac{1}{4}}.
$$
Since $x=\frac{1}{\sqrt{n}}e\ot e$ is rank one, its norm in $S^p$
does not depend on $p$, and it is equal to
$\frac{1}{\sqrt{n}}\norm{e}_{\ell^2_n}^2=\sqrt{n}$. Then, by Lemma
\ref{equivalence usual norm and col for La}, we have
$\norm{x}_{p,\mathcal{L}_a,c,1}\approx  \sqrt{n}$. We obtain the
first equality of (\ref{c/r=infty}) in that case.

If $2<p\leq 4$, we can write
$\frac{1}{\frac{p}{2}}=\frac{1-\theta}{1}+\frac{\theta}{2}$ with
$0<\theta\leq 1$. Then
$$
\norm{x}_{p,\mathcal{L}_a,r,1}^2 =\norm{A}_{S^{\frac{p}{2}}_n} \leq
\norm{A}_{S^1_n}^{1-\theta}\norm{A}_{S^2_n}^{\theta}.
$$
By construction, we have $A\geq 0$, hence we have
\begin{align*}
\norm{A}_{S^1_n}
    &  =\tr\Bigg(\sum_{i,j=1}^n\frac{2^{i+j}}{(2^i+2^j-1)^2}e_{ij}\Bigg)\\
    &  =\sum_{i=1}^n\frac{4^{i}}{(2^{i+1}-1)^2}\leq \sum_{i=1}^n\frac{4^{i}}{(2^{i})^2}= n.
\end{align*}
Thus
$$
\norm{x}_{p,\mathcal{L}_a,r,1}^2 \lesssim
n^{1-\theta}n^{\frac{\theta}{2}}=n^{1-\frac{\theta}{2}}.
$$
Recall that $\norm{x}_{p,\mathcal{L}_a,c,1}\approx \sqrt{n}$. We
obtain that
$$
\frac{\norm{x}_{p,\mathcal{L}_a,c,1}}{\norm{x}_{p,\mathcal{L}_a,r,1}}
\gtrsim\frac{n^{\frac{1}{2}}}{n^{\frac{1}{2}-\frac{\theta}{4}}}=n^{\frac{\theta}{4}}.
$$
Since $n$ was arbitrary and $\theta>0$, we obtain the first part of
(\ref{c/r=infty}) in this case. Likewise, the above proof has a
`right analog' which proves the second equality of
(\ref{c/r=infty}).

We now turn to the proof of (\ref{r/c=infty}). We assume that
$1<p<2$. The second part of (\ref{c/r=infty}) says
\begin{equation}\label{r/c=infty dual}
\sup\Bigg\{\frac{\norm{y}_{p^*,\mathcal{L}_a^*,r,1}}{\norm{y}_{p^*,\mathcal{L}_a^*,c,1}}\
\colon\ y\in S^{p^*}\Bigg\}=\infty.
\end{equation}
To prove the first equality of (\ref{r/c=infty}), assume on the
contrary that there is a constant $K>0$ such that for any $x\in S^p$
\begin{equation}\label{hypothese}
\norm{x}_{p,\mathcal{L}_a,r,1}\leq K\norm{x}_{p,\mathcal{L}_a,c,1}.
\end{equation}

We begin by showing a duality relation between
$\norm{\cdot}_{p^*,\mathcal{L}_a^*,c,1}$ and
$\norm{\cdot}_{p,\mathcal{L}_a,r,1}$. Let $y\in S^{p^*}$ and $x\in
S^{p}$. For any integer $n\geq 1$, recall that $d_n$ is the bounded
diagonal operator on $\ell^2$ defined by the matrix
$\diag(1,\ldots,1,0,\ldots)$. By (\ref{(I-a)^-2=}), for any
$0<\varrho <1$ and any integer $n\geq 1$, we have
\begin{align*}
\MoveEqLeft \left|\big\langle
y,\mathcal{L}_{d_n}(x)\big\rangle_{S^{p^*},S^{p}}\right|
     =\left|\Bigg\langle y,  \sum_{k=1}^{+\infty}k(\varrho \mathcal{L}_a)^{2(k-1)}\big(I-(\varrho \mathcal{L}_a)^2\big)^2\mathcal{L}_{d_n}(x) \Bigg\rangle_{S^{p^*},S^{p}}\right|\\
    &=\left|\sum_{k=1}^{+\infty}\Big\langle y,k(\varrho \mathcal{L}_a)^{2(k-1)}\big(I-(\varrho \mathcal{L}_a)^2\big)^2\mathcal{L}_{d_n}(x) \Big\rangle_{S^{p^*},S^{p}}\right|\\
    &=\left|\sum_{k=1}^{+\infty}\Big\langle k^{\frac{1}{2}}(\varrho \mathcal{L}_a^*)^{k-1}(I-\varrho \mathcal{L}_a^*)(I+\varrho \mathcal{L}_a^*)^2y,k^{\frac{1}{2}}(\varrho \mathcal{L}_a)^{k-1}(I-\varrho \mathcal{L}_a)\mathcal{L}_{d_n}(x) \Big\rangle_{S^{p^*},S^{p}}\right|\\
    &\leq \bgnorm{\Big(k^{\frac{1}{2}}(\varrho \mathcal{L}_a^*)^{k-1}(I-\varrho \mathcal{L}_a^*)(I+\varrho \mathcal{L}_a^*)^2y\Big)_{k\geq 1}}_{S^p(\ell^2_c)}\bnorm{\mathcal{L}_{d_n}(x)}_{p,\varrho\mathcal{L}_a,r,1}.
\end{align*}
Now, it is easy to see that $\{\mathcal{L}_a^*\}$ is Col-bounded. We
infer that
\begin{align*}
\left|\big\langle
y,\mathcal{L}_{d_n}(x)\big\rangle_{S^{p^*},S^{p}}\right|
    &\lesssim \bgnorm{\Big(k^{\frac{1}{2}}(\varrho \mathcal{L}_a^*)^{k-1}\big(I-\varrho \mathcal{L}_a^*)y\Big)_{k\geq 1}}_{S^p(\ell^2_c)}\bnorm{\mathcal{L}_{d_n}(x)}_{p,\varrho \mathcal{L}_a,r,1}\\
    &= \norm{y}_{p^*,\varrho \mathcal{L}_a^*,c,1}\bnorm{\mathcal{L}_{d_n}(x)}_{p,\varrho \mathcal{L}_a,r,1}.
\end{align*}
Assume for a while that $y\in\Ran \big((I-\mathcal{L}_a^*)^2\big)$.
By Lemma \ref{lemma approx} (2), letting $\varrho $ to 1, we obtain
\begin{align*}
\left|\big\langle
y,\mathcal{L}_{d_n}(x)\big\rangle_{S^{p^*},S^{p}}\right|
    &\lesssim
    \norm{y}_{p^*,\mathcal{L}_a^*,c,1}\bnorm{\mathcal{L}_{d_n}(x)}_{p,\mathcal{L}_a,r,1}.
\end{align*}
Letting $n$ to the infinity, we obtain
\begin{equation*}
\left|\langle y,x\rangle_{S^{p^*},S^{p}}\right|\lesssim
\norm{y}_{p^*,\mathcal{L}_a^*,c,1}\norm{x}_{p,\mathcal{L}_a,r,1}.
\end{equation*}

According to (\ref{hypothese}) and the first part of
(\ref{equivalence column square function}), we deduce that
\begin{align*}
 \left|\langle y,x \rangle_{S^{p^*},S^{p}}\right|
    & \lesssim \norm{y}_{p^*,\mathcal{L}_a^*,c,1}\norm{x}_{p,\mathcal{L}_a,c,1}\\
    & \lesssim \norm{y}_{p^*,\mathcal{L}_a^*,c,1}\norm{x}_{S^p}.
\end{align*}
By duality, we finally obtain that
\begin{equation}
\label{inequality 4.9} \norm{y}_{S^{p^*}}\lesssim
\norm{y}_{p^*,\mathcal{L}_a^*,c,1}.
\end{equation} For an arbitrary
$y\in S^{p^*}$, we also obtain (\ref{inequality 4.9}) by applying it
to $\mathcal{L}_{d_n}^*(y)$ and then passing to the limit. The
second equivalence of (\ref{equivalence column square function})
says that $\norm{y}_{p^*,\mathcal{L}_a^*,r,1}\approx
\norm{y}_{S^{p^*}}$ for any $y\in S^{p^*}$. This contradicts
(\ref{r/c=infty dual}) and completes the proof of the first part of
(\ref{r/c=infty}). The proof of the second part is similar.
\end{preuve}

For a operator admitting a completely bounded $H^\infty(B_\gamma)$
functional calculus, it also seems interesting, in view of the
equivalence (\ref{equivalence usual norm and T1}), to compare the
column and row square functions with the usual norm
$\norm{\cdot}_{L^p(M)}$. If $T$ is a operator with $\Ran(I-T)$ dense
in $L^p(M)$ which admits a bounded $H^\infty(B_\gamma)$ functional
calculus for some $\gamma\in \big]0,\frac{\pi}{2}\big[$, the
equivalence (\ref{equivalence usual norm and T1}) and Theorems
\ref{principal thm LM ghost} and \ref{thm principal ALM} implies
that
$$
 \norm{x}_{L^p(M)}\lesssim \norm{x}_{p,T,c,1}\ \ \ \ \ \text{and} \ \
\ \ \ \norm{x}_{L^p(M)}\lesssim \norm{x}_{p,T,r,1}
$$
if $1 < p \leq  2$ and
$$
\norm{x}_{p,T,c,1} \lesssim \norm{x}_{L^p(M)}\ \ \ \ \ \text{and} \
\ \ \ \ \norm{x}_{p,T,r,1} \lesssim \norm{x}_{L^p(M)}
$$
if $2 \leq p<\infty $, for any $x\in L^p(M)$. The following result
says that except for $p=2$, these estimates cannot be reversed:
\begin{cor}
Suppose that $2<p<\infty$ (resp. $1<p<2$). Let $\alpha> 0$. There
exists a Ritt operator $T$ on the Schatten space $S^p$, with
$\Ran(I-T)$ dense in $S^p$, which admits a completely bounded
$H^\infty(B_\gamma)$ functional calculus with $\gamma\in
\big]0,\frac{\pi}{2}\big[$ such that
\begin{equation*}
\sup\bigg\{\frac{\norm{x}_{S^p}}{\norm{x}_{p,T,c,\alpha}}\ \colon\
x\in S^p\bigg\}=\infty\ \ \Bigg(\text{resp.
}\sup\bigg\{\frac{\norm{x}_{p,T,c,\alpha}}{\norm{x}_{S^p}}\ \colon\
x\in S^p\bigg\}=\infty\Bigg).
\end{equation*}
Moreover, the same result holds with $\norm{\cdot}_{p,T,c,\alpha}$
replaced by $\norm{\cdot}_{p,T,r,\alpha}$.
\end{cor}

\begin{preuve}
One more time, we only need to prove this result for $\alpha=1$.
Then, this follows from Lemma \ref{equivalence usual norm and col
for La} and Theorem \ref{thm non equivalence}.
\end{preuve}


\section{An alternative square function for $1<p<2$}


Let $T$ be a Ritt operator on $L^p(M)$, with $1<p<2$. For any
$\alpha>0$, we may consider an alternative square function by
letting
$$
\norm{x}_{p,T,0,\alpha}=\inf\Big\{\norm{x_1}_{p,T,c,\alpha}+\norm{x_2}_{p,T,r,\alpha}\
\colon\ x=x_1+x_2 \Big\}
$$
for any $x\in L^p(M)$.

Note that if $T$ is both Col-Ritt and Row-Ritt, by Theorem \ref{Thm
equivalence square function colrow}, the square functions
$\norm{x}_{p,T,0,\alpha}$ and $\norm{x}_{p,T,0,\beta}$ are
equivalent for any $\alpha,\beta>0$.

Suppose that $\norm{x}_{p,T,0,\alpha}$ is finite and that we have a
decomposition $x=x_1+x_2$ with $\norm{x_1}_{p,T,c,\alpha} <\infty$
and $\norm{x_2}_{p,T,r,\alpha}<\infty$. Letting
$u_k=k^{\alpha-\frac{1}{2}}T^{k-1}(I-T)^{\alpha}x_1$ and
$v_k=k^{\alpha-\frac{1}{2}}T^{k-1}(I-T)^{\alpha}x_2$, we have
$$
k^{\alpha-\frac{1}{2}}T^{k-1}(I-T)^{\alpha}x=u_k+v_k, \ \ \ \ \ \ k
\geq 1.
$$
Moreover, the sequences $u$ and $v$ belong to
$L^p\big(M,\ell^2_c\big)$ and $L^p\big(M,\ell^2_r\big)$
respectively. We deduce that
$$
\norm{x}_{p,T,\alpha} \leq \norm{x}_{p,T,0,\alpha}, \ \ \ \ \ x\in
L^p(M).
$$
We do not know if the two square functions
$\norm{\cdot}_{p,T,\alpha}$ and $\norm{\cdot}_{p,T,0,\alpha}$ are
equivalent in general. In the next statement, we give a sufficient
condition for an such equivalence to hold true.

\begin{thm}
\label{equivalence norm and norm0} Suppose $1<p<2$. Let $T$ be a
bounded operator on $L^p(M)$ with $\Ran(I-T)$ dense in $L^p(M)$.
Assume that $T$ is both Col-Ritt and Row-Ritt. Let $\alpha,\eta>0$.
Suppose that $T$ satisfies a `dual square function estimate'
\begin{equation}\label{dual square estimate}
\norm{y}_{p^*,T^*,\eta}\lesssim \norm{y}_{L^{p^*}(M)},\qquad y\in
L^{p^*}(M).
\end{equation}
Then we have an equivalence
$$
\norm{x}_{p,T,\alpha}\approx \norm{x}_{p,T,0,\alpha},\qquad x\in
L^p(M).
$$
Indeed, there is a positive constant $C$ such that whenever $x\in
L^p(M)$ satisfies $\norm{x}_{p,T,\alpha}<\infty$, then there exists
$x_1,x_2\in L^p(M)$ such that
$$
x=x_1+x_2\ \ \text{and} \ \
\norm{x_1}_{p,T,c,\alpha}+\norm{x_2}_{p,T,r,\alpha}\leq C
\norm{x}_{p,T,\alpha}.
$$
\end{thm}

\begin{preuve}
Since $T$ is both Col-Ritt and Row-Ritt, it is also an $R$-Ritt
operator. Then, by Theorem \ref{thm principal ALM} and Theorem
\ref{Thm equivalence square function colrow}, we only need to prove
this result for $\alpha=1$ and $\eta=1$.
Observe that, for any $y\in L^{p^*}(M)$, we have
\begin{align*}
\MoveEqLeft
\bgnorm{\Big(k^{\frac{1}{2}}(T^*)^{k-1}(I+T^*)^{2}(I-T^*)^{}y\Big)_{k\geq 1}}_{L^{p^*}(M,\ell^2_{\rad})}\\
    &\lesssim \bnorm{(I+T^*)^{2}}_{L^{p^*}(M)\to L^{p^*}(M)} \bgnorm{\Big(k^{\frac{1}{2}}(T^*)^{k-1}(I-T^*) y\Big)_{k\geq 1}}_{L^{p^*}(M,\ell^2_{\rad})}\\
    &\lesssim \norm{y}_{L^{p^*}(M)}\hspace{1cm} \text{by (\ref{dual square estimate})}.
\end{align*}
We let
$$
\begin{array}{cccc}
 Z:   &  L^{p^*}(M)   &  \longrightarrow   &  L^{p^*}\big(M,\ell^2_{\rad}\big)  \\
      &      y        &  \longmapsto       & \displaystyle \Big(k^{\frac{1}{2}}(T^*)^{k-1}(I+T^*)^{2}(I-T^*)^{}y\Big)_{k\geq 1} \\
\end{array}
$$
denote the resulting bounded map. Let $x\in L^p(M)$ such that
$\norm{x}_{p,T,1}<\infty$. There exists two elements $u\in
L^p\big(M,\ell^2_c\big)$ and $v \in L^p\big(M,\ell^2_r\big)$ such
that for any positive integer $k$
\begin{equation}\label{uk+vk}
    u_k+v_k=k^{\frac{1}{2}}T^{k-1}(I-T)^{}x
\end{equation}
and such that
$$
\norm{u}_{L^p(M,\ell^2_c)}+\norm{v}_{L^p(M,\ell^2_r)}\leq
2\norm{x}_{p,T,1}.
$$
Recall that we have contractive inclusions
$L^{p}\big(M,\ell^2_{c}\big) \subset L^{p}\big(M,\ell^2_{\rad}\big)$
and $L^{p}\big(M,\ell^2_{r}\big) \subset
L^{p}\big(M,\ell^2_{\rad}\big)$. Thus, by (\ref{dualite
Lp(M,ell2rad)}), we can define $x_1$ and $x_2$ of $L^p(M)$ by
$$
x_1=Z^*u \ \ \ \ \ \text{and} \ \ \ \ \ x_2=Z^*v.
$$
We will show that $x=x_1+x_2$. Since $T$ is a Col-Ritt-operator, by
Proposition \ref{prop bornitude} (or by \cite[Proposition
2.8]{ALM}), we infer that there exists a positive constant $C$ such
that
\begin{align*}
 \sum_{k=1}^{+\infty} \Bnorm{k^{\frac{1}{2}}T^{k-1}(I-T)^{2}}_{L^p(M) \to L^p(M)}^2
   &=  \sum_{k=1}^{+\infty}k^{}\Bnorm{T^{k-1}(I-T)^{2}}_{L^p(M) \to L^p(M)}^2  \\
   &\leq C^2\sum_{k=1}^{+\infty}\frac{1}{k^{3}}< \infty.
\end{align*}
For any $1 < p < 2$, by \cite[Proposition 2.5]{JMX}, we have the
contractive inclusion $L^{p}\big(M,\ell^2_{c}\big) \subset
\ell^2\big(L^p(M)\big)$. We deduce that $\sum_{k=1}^{+\infty}
\norm{u_k}_{L^p(M)}^2 < \infty.$ According to the Cauchy-Schwarz
inequality, we deduce that the series
$$
\sum_{k=1}^{+\infty}k^{\frac{1}{2}}T^{k-1}\big(I-T^2\big)^{2}u_k =
(I+T)^{2}\sum_{k=1}^{+\infty}k^{\frac{1}{2}}T^{k-1}(I-T)^{2}u_k
$$
converges absolutely in $L^p(M)$. Now, for any $y\in L^{p^*}(M)$, we
have
\begin{align*}
\Big\langle(I-T)x_1,y\Big\rangle_{L^{p}(M),L^{p^*}(M)}
   &= \Big\langle(I-T)^{}Z^*u,y\Big\rangle_{L^{p}(M),L^{p^*}(M)}\\
   &= \Big\langle u, Z(I-T^{*})^{}y \Big\rangle_{L^{p}(M,\ell^2_\rad),L^{p^*}(M,\ell^2_\rad)} \\
   &= \bigg\langle u,\Big(k^{\frac{1}{2}}(T^*)^{k-1} (I+T^*)^{2}(I-T^*)^{2}y\Big)_{k\geq 1}\bigg\rangle_{L^{p}(M,\ell^2_\rad),L^{p^*}(M,\ell^2_\rad)}\\
   &= \sum_{k=1}^{+\infty}\Big\langle  u_k,k^{\frac{1}{2}}(T^*)^{k-1} \big(I-(T^*)^2\big)^{2}y\Big\rangle_{L^{p}(M),L^{p^*}(M)}\\
   &= \Bigg\langle\sum_{k=1}^{+\infty} k^{\frac{1}{2}}T^{k-1} \big(I-T^2\big)^{2}u_k ,y\Bigg\rangle_{L^{p}(M),L^{p^*}(M)}.
\end{align*}
Thus, we deduce that
\begin{equation}\label{(I-T)x1}
(I-T)^{}x_1=\sum_{k=1}^{+\infty}
k^{\frac{1}{2}}T^{k-1}\big(I-T^2\big)^{2}u_k.
\end{equation}
Similarly we have
\begin{equation*}\label{(I-T)x2}
(I-T)^{}x_2=\sum_{k=1}^{+\infty}
k^{\frac{1}{2}}T^{k-1}\big(I-T^2\big)^{2}v_k.
\end{equation*}
Now, we infer that
\begin{align*}
(I-T)^{}(x_1+x_2)
   &= \sum_{k=1}^{+\infty}  k^{\frac{1}{2}}T^{k-1}\big(I-T^2\big)^{2}u_k +\sum_{k=1}^{+\infty} k^{\frac{1}{2}}T^{k-1}\big(I-T^2\big)^{2}v_k \\
   &= \sum_{k=1}^{+\infty}  k^{\frac{1}{2}}T^{k-1}\big(I-T^2\big)^{2}(u_k+v_k)\\
   &= \sum_{k=1}^{+\infty}  k^{\frac{1}{2}}T^{k-1}\big(I-T^2\big)^{2}k^{\frac{1}{2}}T^{k-1}(I-T)^{}x\hspace{1cm} \text{by (\ref{uk+vk})}\\
   &= \sum_{k=1}^{+\infty}  k T^{2k-2}(I+T)^{2}(I-T)^{3}x.
\end{align*}
By (\ref{power series}), for any $z\in \mathbb{D}$, we have
$$
\sum_{k=1}^{+\infty}kz^{2k-2}(1-z^2)^2 = 1.
$$
Since the operator $T$ is power bounded, we note that for every
$\varrho\in ]0,1[$ we have
\begin{equation}\label{formule pour Id}
I=\sum_{k=1}^{+\infty}k(\varrho T)^{2k-2}\big(I-(\varrho
T)^2\big)^{2},
\end{equation}
the series being absolutely convergent. Hence, for any $\varrho\in
]0,1[$, we have
\begin{align*}
(I-\varrho T)x
   &=(I-\varrho T)\sum_{k=1}^{+\infty}k(\varrho T)^{2k-2}\big(I-(\varrho T)^2\big)^{2}x \\
   &= \sum_{k=1}^{+\infty}k(\varrho T)^{2k-2}(I+\varrho T)^{2}(I-\varrho T)^{3}x.
\end{align*}
It is not difficult to see that the latter series is normally
convergent on [0,1]. Hence, letting $\varrho$ to $1$, we deduce that
$$
(I-T)x= \sum_{k=1}^{+\infty}kT^{2k-2}(I+T)^{2}(I-T)^{3}x.
$$
Then we obtain
$$
(I-T)x=(I-T)(x_1+x_2).
$$
Since the space $\Ran(I-T)$ is dense in $L^p(M)$, by the Mean
Ergodic Theorem (see \cite[Section 2.1]{Kre}), the operator $I-T$ is
injective. Consequently, we have $x=x_1+x_2$. Now, it remains to
estimate $\norm{x_1}_{p,T,1,c}$ and $\norm{x_2}_{p,T,1,r}$.
According to (\ref{(I-T)x1}), we have
\begin{align*}
m^{\frac{1}{2}}T^{m-1}(I-T) x_1 &= \sum_{k=1}^{+\infty}
k^{\frac{1}{2}} m^{\frac{1}{2}}T^{k+m-2}\big(I-T^2\big)^{2}u_k
\end{align*}
for any integer $m\geq 1$. It is convenient to write this as
$m^{\frac{1}{2}}T^{m-1}(I-T) x_1=(I+T)^2y_{m} $ with
\begin{equation}
\label{y_m} y_{m}=\sum_{k=1}^{+\infty}
k^{\frac{1}{2}}m^{\frac{1}{2}}T^{k+m-2} (I-T)^{2}u_k.
\end{equation}
Now, observe that
\begin{align*}
k^{\frac{1}{2}} m^{\frac{1}{2}}T^{k+m-2} (I-T)^{2}
&=\frac{k^{\frac{1}{2}} m^{\frac{1}{2}}}{(k+m-1)^{2}}\cdot
(k+m-1)^{2}T^{k+m-2}(I-T)^{2}.
\end{align*}
According to \cite[Proposition 2.3]{ALM} and \cite[Lemma 2.4]{ALM},
the matrix
$$
\biggl[\frac{k^{\frac{1}{2}}
m^{\frac{1}{2}}}{(k+m-1)^{2}}\biggr]_{k,m\geq 1}
$$
represents an element of $B(\ell^2)$. Moreover, by Proposition
\ref{prop bornitude}, the set
$$
\Big\{ (k+m-1)^{2}T^{k+m-2}(I-T)^{2}\,:\, k,m\geq 1 \Big\}
$$
is Col-bounded. By Proposition \ref{prop 2.5 version col}, we deduce
that $(y_m)_{m\geq 1} \in L^p\big(M,\ell^2_c\big)$ and that
\begin{align*}
\bnorm{(y_{m})_{m \geq 1}}_{L^p(M,\ell^2_c)}
&\lesssim \norm{u}_{L^p(M,\ell^2_c)}.
\end{align*}
Since $\{T\}$ is Col-bounded, we have
\begin{align*}
\norm{x_1}_{p,T,c,1}
   &= \bgnorm{\Big(m^{\frac{1}{2}}T^{m-1}(I-T)^{}x_1\Big)_{m\geq 1}}_{L^p(M,\ell^2_c)}\\
   &= \bgnorm{\Big((I+T)^{2}y_m\Big)_{m\geq 1}}_{L^p(M,\ell^2_c)} \hspace{1cm} \text{by (\ref{y_m})} \\
   &\lesssim  \bnorm{(y_{m})_{m \geq 1}}_{L^p(M,\ell^2_c)}.
\end{align*}
Finally, we deduce that there exists a positive constant $C$ such
that
\begin{align*}
\norm{x_1}_{p,T,c,1}  &\leq C \norm{u}_{L^p(M,\ell^2_c)}.
\end{align*}
Moreover, we have a similar result for $x_2$. Finally, we have
\begin{align*}
\norm{x_1}_{p,T,c,1}+\norm{x_2}_{p,T,r,1}
    &\leq C\norm{u}_{L^p(M,\ell^2_c)}+C\norm{v}_{L^p(M,\ell^2_r)}\\
    &\leq C\norm{x}_{p,T,1}.
\end{align*}
\end{preuve}

\begin{cor}
\label{coro CB funct calculus implies equivalence} Suppose $1<p<2$.
Let $T$ be a bounded operator on $L^p(M)$ with $\Ran(I-T)$ dense in
$L^p(M)$ and let $\alpha>0$. Assume that $T$ admits a completely
bounded $H^\infty(B_\gamma)$ functional calculus for some $\gamma
\in \big]0,\frac{\pi}{2}\big[$. Then we have an equivalence
$$
\inf\Big\{\norm{x_1}_{p,T,c,\alpha}+\norm{x_2}_{p,T,r,\alpha}\
\colon\ x=x_1+x_2 \Big\}\approx\norm{x}_{L^p(M)},\qquad x \in
L^p(M).
$$
\end{cor}

\begin{preuve}
By Theorem \ref{FC CB implies Col}, the operator $T$ is both
Col-Ritt and Row-Ritt (hence $R$-Ritt). Moreover, by Theorem
\ref{principal thm LM ghost}, it satisfies a `dual square estimate'
$$
\norm{y}_{p^*,T^*,1}\lesssim \norm{y}_{L^{p^*}(M)},\qquad y\in
L^{p^*}(M).
$$
Then, by Theorem \ref{equivalence norm and norm0} above, the norms
$\norm{\cdot}_{p,T,\alpha}$ and $\norm{\cdot}_{p,T,0,\alpha}$ are
equivalent. Furthermore, by Theorem \ref{thm principal ALM} and
(\ref{equivalence usual norm and T1}), $\norm{\cdot}_{p,T,\alpha}$
is equivalent to the usual norm $\norm{\cdot}_{L^{p}(M)}$, which
proves the result.
\end{preuve}

Assume now that $\tau$ is finite and normalized, that is, $\tau(1) =
1$. Following \cite{HaM} and \cite{Ric} (see also \cite{AD}), we say
that a linear map $T$ on $M$ is a Markov map if T is unital,
completely positive and trace preserving. As is well known, such a
map is necessarily normal and for any $1 \leq p < \infty$, it
extends to a contraction $T_p$ on $L^p(M)$. We say that $T$ is
selfadjoint if, for any $x,x'\in M$, we have
$$
\tau\big(T(x)x'\big)=\tau\big(xT(x')\big).
$$
This is equivalent to $T_2$ being selfadjoint in the Hilbertian
sense. We also consider the operator
$$
A_p = I_{}- T_p.
$$
The following result is proved in the proof of \cite[Proposition
8.7]{LM0} with \textit{bounded} instead of \textit{completely
bounded}. But a careful reading of the proof shows that we have this
stronger result. We refer to \cite{Haa}, \cite{JMX}, \cite{LM3} and
\cite{LM0} for information on $H^\infty(\Sigma_\theta)$ functional
calculus.

\begin{prop}
\label{prop Markov implies CB} Suppose $1 < p < \infty$. Let $T$ be
a selfadjoint Markov map on $M$. Then the operator $A_p$ is
sectorial and admits a completely bounded $H^\infty(\Sigma_\theta)$
functional calculus for some $\theta \in \big]
0,\frac{\pi}{2}\big[$.
\end{prop}

Assume $1<p<\infty$. At this point, it is crucial to recall that
$L^p$-realizations $T_p$ of Markov maps $T$ on $M$ such that
$-1\notin \sigma(T_2)$ are Ritt operators, as noticed by C. Le Merdy
in \cite{LM0}. Let $T$ be a selfadjoint Markov map on $M$. According
to \cite{LM0} and Proposition \ref{prop Markov implies CB}, we
obtain that $T_p$ admits a completely bounded $H^\infty(B_\gamma)$
functional calculus for some $\gamma \in \big]0,\frac{\pi}{2}\big[$.
Hence, by Corollary \ref{coro CB funct calculus implies
equivalence}, we deduce the following result which strengthens a
result of \cite{LM0}.


\begin{cor}
Suppose $1<p< 2$. Let $T$ be a selfadjoint Markov map on $M$ such
that $-1\notin \sigma(T_2)$ with $\Ran(I-T_p)$ dense in $L^p(M)$.
Then, for any $\alpha> 0$ there exists a positive constant $C$ such
that for any $x \in L^p(M)$, there exists $x_1, x_2 \in L^p(M)$
satisfying $x=x_1+x_2$ and
$$
\Bgnorm{\bigg(\sum_{k=1}^{+\infty}k^{2\alpha-1}\left|T^{k-1}(I-T)^{\alpha}(x_1)\right|^2\bigg)^{\frac{1}{2}}}_{L^p}
+\Bgnorm{\bigg(\sum_{k=1}^{+\infty}k^{2\alpha-1}
\left|\Big(T^{k-1}(I-T)^{\alpha}(x_2)\Big)^*\right|^2\bigg)^{\frac{1}{2}}}_{L^p}
\leq C \norm{x}_{L^p(M)}.
$$
\end{cor}

\bigskip


\textbf{Acknowledgement}. We wish to thank my thesis adviser
Christian Le Merdy for his support and advice and \'Eric Ricard for
simplifications in some approximation arguments.


\bigskip\footnotesize{
\n Laboratoire de Math\'ematiques, Universit\'e de Franche-Comt\'e,
25030 Besan\c{c}on Cedex,  France\\
cedric.arhancet@univ-fcomte.fr\hskip.3cm

\end{document}
For a Ritt operator $T$, we let $P_T$ denote the projection onto
$\Ker (I-T)$ which vanishes on $\ovl{\Ran(I-T)}$.

\begin{prop}
\label{equivalence norm projection} Let $T$ be a Ritt operator on
$L^p(M)$. Assume that $T$ admits a completely bounded
$H^\infty(B_\gamma)$ functional calculus for some $\gamma\in
\big]0,\frac{\pi}{2}\big[$. Then we have an equivalence
$$
\norm{x}_{L^p(M)} \approx \bnorm{P_T(x)}_{L^p(M)}+ \norm{x}_{p,T,1},
\qquad x \in L^p(M).
$$
\end{prop}

\begin{proof}
Comme le cas commutatif de Le Merdy Xu 1.
\end{proof}
we have CHANGER PRESENTATION?
\begin{align*}
\norm{x}_{\varrho\mathcal{L}_a,c,1}
&=\Bgnorm{\bigg(\sum_{k=1}^{+\infty} k\big((\varrho \mathcal{L}_{a})^{k-1}(I-\varrho\mathcal{L}_{a})(x)\big)^{*}\big((\varrho \mathcal{L}_{a})^{k-1}(I-\varrho \mathcal{L}_{a})(x)\big)\bigg)^{\frac{1}{2}}}_{S^p} \\
&=\Bgnorm{\bigg(\sum_{k=1}^{+\infty} k\big((\varrho a)^{k-1}(I-\varrho a)x\big)^{*}\big((\varrho a)^{k-1}(I-\varrho a)x\big)\bigg)^{\frac{1}{2}}}_{S^p}\\
&=\Bgnorm{\bigg(\sum_{k=1}^{+\infty} kx^*(I-\varrho a)(\varrho
a)^{2(k-1)}(I-\varrho a)x\bigg)^\frac{1}{2}}_{S^p}.
\end{align*}
Now, for any $z\in \mathbb{D}$, we have
\begin{equation}
\label{power series}
    \sum_{k=1}^{+\infty} ka^{k-1}=(1-a)^{-2}.
\end{equation}
Since the operator $a$ is a contraction, we deduce that, for every
$\varrho \in ]0,1[$, the operator $I-(\varrho a)^2$ is inversible
and that we have
\begin{equation}\label{(I-a)^-2=}
   \sum_{k=1}^{+\infty}k (\varrho a)^{2(k-1)}=\big(I-(\varrho a)^2\big)^{-2}
\end{equation}
the series being absolutely convergent. Then we deduce that
\begin{align*}
\norm{x}_{\varrho\mathcal{L}_a,c,1} &=\Bgnorm{\Bigg(x^*(I-\varrho
a)^2\bigg(\sum_{k=1}^{+\infty}k(\varrho
a)^{2(k-1)}\bigg)x\Bigg)^{\frac{1}{2}}}_{S^p}\\
&=\bgnorm{\Big(x^*(I-\varrho a)^2\big(I-(\varrho a)^2\big)^{-2}x\Big)^{\frac{1}{2}}}_{S^p}\\
&=\Bnorm{\big(x^*(I-\varrho a)^2(I-\varrho a)^{-2}(I+\varrho a)^{-2}x\big)^{\frac{1}{2}}}_{S^p}\\
&=\Bnorm{\big(x^*(I+\varrho a)^{-2}x\big)^{\frac{1}{2}}}_{S^p}\\
&=\Bnorm{\big((I+\varrho a)^{-1}x\big)^*(I+\varrho a)^{-1}x\big)^{\frac{1}{2}}}_{S^p}\\
&=\bnorm{(I+\varrho a)^{-1}x}_{S^p}.
\end{align*}
\begin{align*}
\norm{x}_{p,\mathcal{L}_a,r,1}
&=\Bgnorm{\bigg(\sum_{k=1}^{+\infty} k\big(\mathcal{L}_{a}^{k-1}(I-\mathcal{L}_{a})(x)\big)\big(\mathcal{L}_{a}^{k-1}(I-\mathcal{L}_{a})(x)\big)^{*}\bigg)^{\frac{1}{2}}}_{S^p} \\
&=\Bgnorm{\bigg(\sum_{k=1}^{+\infty} k\big(a^{k-1}(I-a)x\big)\big(a^{k-1}(I-a)x\big)^{*}\bigg)^{\frac{1}{2}}}_{S^p}\\
&=\Bgnorm{\bigg(\sum_{k=1}^{+\infty} ka^{k-1}(I-a)xx^*(I-a)a^{k-1}\bigg)^\frac{1}{2}}_{S^p}\\
&=\Bgnorm{\bigg(\sum_{i,j=1}^n\sum_{k=1}^{+\infty} ka^{k-1}(I-a)e_{ij}(I-a)a^{k-1}\bigg)^\frac{1}{2}}_{S^p}\\
&=\Bgnorm{\bigg(\sum_{i,j=1}^n\sum_{k=1}^{+\infty} ka_i^{k-1}(1-a_i)(1-a_j)a_j^{k-1}e_{ij}\bigg)^\frac{1}{2}}_{S^p}\\
&=\Bgnorm{\bigg(\sum_{i,j=1}^n(1-a_i)(1-a_j)\sum_{k=1}^{+\infty}
k(a_ia_j)^{k-1}e_{ij}\bigg)^\frac{1}{2}}_{S^p}.
\end{align*}
Using the equality (\ref{power series}), we obtain
\begin{align*}
\norm{x}_{p,\mathcal{L}_a,r,1}
&=\Bgnorm{\bigg(\sum_{i,j=1}^n(1-a_i)(1-a_j)(1-a_ia_j)^{-2}e_{ij}\bigg)^\frac{1}{2}}_{S^p}\\
&=\Bgnorm{\bigg(\sum_{i,j=1}^n\frac{1}{2^{i+j}\Big(\frac{1}{2^i}+\frac{1}{2^j}-\frac{1}{2^{i+j}}\Big)^2}e_{ij}\bigg)^\frac{1}{2}}_{S^p}\\
&=\Bgnorm{\bigg(\sum_{i,j=1}^n\frac{2^{i+j}}{(2^i+2^j-1)^2}e_{ij}\bigg)^\frac{1}{2}}_{S^p}\\
&=\Bgnorm{\sum_{i,j=1}^n\frac{2^{i+j}}{(2^i+2^j-1)^2}e_{ij}}_{S^{\frac{p}{2}}}^\frac{1}{2}.
\end{align*}

\section{Alternative square function for $1<p<2$.}


Let $T$ a Ritt operator on $L^p(M)$, with $1<p<2$. For any
$\alpha>0$, we may consider an alternative square function by
letting
$$
\norm{x}_{p,T,0,\alpha}=\inf\Big\{\norm{x_1}_{p,T,c,\alpha}+\norm{x_2}_{p,T,r,\alpha}\
\colon\ x=x_1+x_2 \Big\}
$$
for any $x\in L^p(M)$. Suppose that $\norm{x}_{p,T,0,\alpha}$ is
finite and that we have a decomposition $x=x_1+x_2$ with
$\norm{x_1}_{p,T,c,\alpha} <\infty$ and
$\norm{x_2}_{p,T,r,\alpha}<\infty$. Letting
$u_k=k^{\alpha-\frac{1}{2}}T^{k-1}(I-T)^{\alpha}x_1$ and
$v_k=k^{\alpha-\frac{1}{2}}T^{k-1}(I-T)^{\alpha}x_2$, we have
$$
k^{\alpha-\frac{1}{2}}T^{k-1}(I-T)^{\alpha}x=u_k+v_k, \ \ \ \ \ \ k
\geq 1.
$$
Moreover, the sequences $u$ and $v$ belong to
$L^p\big(M,\ell^2_c\big)$ and $L^p\big(M,\ell^2_r\big)$
respectively. We deduce that
$$
\norm{x}_{p,T,\alpha} \leq \norm{x}_{p,T,0,\alpha}, \ \ \ \ \ x\in
L^p(M).
$$
We do not know if the two square functions
$\norm{\cdot}_{p,T,\alpha}$ and $\norm{\cdot}_{p,T,0,\alpha}$ are
equivalent in general. In the next statement, we give a sufficient
condition for an such equivalence to hold true.

\begin{thm}
\label{equivalence norm and norm0} Suppose $1<p<2$. Let $T$ be an
operator on $L^p(M)$ with $R(I-T)$ dense in $L^p(M)$. Assume that
$T$ is both Col-Ritt and Row-Ritt. Let $\alpha,\eta>0$. Suppose that
$T$ satisfies a `dual square function estimate'
\begin{equation}\label{dual square estimate}
\norm{y}_{p^*,T^*,\eta}\lesssim \norm{y}_{L^{p^*}(M)},\qquad y\in
L^{p^*}(M).
\end{equation}
Then we have an equivalence
$$
\norm{x}_{p,T,\alpha}\approx \norm{x}_{p,T,0,\alpha},\qquad x\in
L^p(M).
$$
Indeed, there is a constant $C\geq 1$ such that whenever $x\in
L^p(M)$ satisfies $\norm{x}_{p,T,\alpha}<\infty$, then there exists
$x_1,x_2\in L^p(M)$ such that
$$
x=x_1+x_2\ \ \text{and} \ \
\norm{x_1}_{p,T,c,\alpha}+\norm{x_2}_{p,T,r,\alpha}\leq C
\norm{x}_{p,T,\alpha}.
$$
\end{thm}

\begin{proof}
Let $x\in L^p(M)$ such that $\norm{x}_{p,T,\alpha}<\infty$. There
exists two elements $u\in L^p\big(M,\ell^2_c\big)$ and $v \in
L^p\big(M,\ell^2_r\big)$ such that for any positive integer $k$
\begin{equation}\label{uk+vk}
    u_k+v_k=k^{\alpha-\frac{1}{2}}T^{k-1}(I-T)^{\alpha}x
\end{equation}
and such that
$$
\norm{u}_{L^p(M,\ell^2_c)}+\norm{v}_{L^p(M,\ell^2_r)}\leq
2\norm{x}_{p,T,\alpha}.
$$
By Theorem \ref{thm principal ALM}, the square functions
$\norm{\cdot}_{p^*,T^*,\eta}$ and $\norm{\cdot}_{T^*,\beta}$ are
equivalent for any $\eta,\beta>0$. Then we can suppose $\eta$ such
that $\eta+\alpha$ is an integer $N\geq 1$.

Now, for any integer $k\geq 1$, we define the complex number
$$
c_k=\frac{k(k+1)\cdots (k+N-2)}{k^{\alpha-\frac{1}{2}}}
$$
with the convention that $c_k=\frac{1}{k^{\alpha-\frac{1}{2}}}$ if
$N=1$. We
have $c_k\sim_{+\infty} k^{\eta-\frac{1}{2}}$. 
It follows that, for any $y\in L^{p^*}(M)$, we have
\begin{align*}
\MoveEqLeft
\Bgnorm{\Bigg(\frac{c_k(T^*)^{k-1}(I-(T^*)^2)^{\eta}(I+T^*)^{\alpha}y}{(N-1)!}\Bigg)_{k\geq 1}}_{L^{p^*}(M,\ell^2_{\rad})}\\
    &\lesssim \bgnorm{\frac{1}{(N-1)!}(I+T^*)^{\eta+\alpha}}_{L^{p^*}(M)\to L^{p^*}(M)} \bgnorm{\Big(k^{\eta-\frac{1}{2}}(T^*)^{k-1}(I-T^*)^\eta y\Big)_{k\geq 1}}_{L^{p^*}(M,\ell^2_{\rad})}\\
    &\lesssim \norm{y}_{L^{p^*}(M)}\hspace{1cm} \text{by (\ref{dual square estimate})}.
\end{align*}
Then, we obtain the bounded linear map
$$
\begin{array}{cccc}
 Z:   &  L^{p^*}(M)   &  \longrightarrow   &  L^{p^*}\big(M,\ell^2_{\rad}\big)  \\
      &      y        &  \longmapsto       & \displaystyle \Bigg(\frac{c_k(T^*)^{k-1}(I-(T^*)^2)^{\eta}(I+T^*)^{\alpha}y}{(N-1)!}\Bigg)_{k\geq 1}.  \\
\end{array}
$$
Recall that we have contractive inclusions
$L^{p}\big(M,\ell^2_{c}\big) \subset L^{p}\big(M,\ell^2_{\rad}\big)$
and $L^{p}\big(M,\ell^2_{r}\big) \subset
L^{p}\big(M,\ell^2_{\rad}\big)$. We define the elements $x_1$ and
$x_2$ of $L^p(M)$ by
$$
x_1=Z^*u \ \ \ \ \ \text{and} \ \ \ \ \ x_2=Z^*v.
$$
We will show that $x=x_1+x_2$. Since $T$ is a Ritt-operator, there
exists a positive constant $C$ such that
\begin{align*}
 \sum_{k=1}^{+\infty} \Bnorm{k^{\eta-\frac{1}{2}}T^{k-1}(I-T)^{\eta+1}}_{L^p(M) \to L^p(M)}^2
   &=  \sum_{k=1}^{+\infty}k^{2\eta-1}\bnorm{T^{k-1}(I-T)^{\eta+1}}_{L^p(M) \to L^p(M)}^2  \\
   &\leq  C^2\sum_{k=1}^{+\infty} k^{2\eta-1}\bigg(\frac{1}{k^{\eta+1}}\bigg)^{2}  \\
   &=C^2\sum_{k=1}^{+\infty}\frac{1}{k^{3}}< \infty.
\end{align*}
By \cite[Proposition 2.5]{JMX}, we have the contractive inclusion
$L^{p}\big(M,\ell^2_{c}\big) \subset \ell^2\big(L^p(M)\big)$. We
deduce that $\sum_{k=1}^{+\infty} \norm{u_k}_{L^p(M)}^2 < \infty.$
According to Cauchy-Schwarz inequality, we deduce that the series
$$
\sum_{k=1}^{+\infty}k^{\eta-\frac{1}{2}}T^{k-1}(I-T^2)^{\eta}(I+T)^{\alpha}(I-T)^{}u_k
=
(I+T)^{\eta+\alpha}\sum_{k=1}^{+\infty}k^{\eta-\frac{1}{2}}T^{k-1}(I-T)^{\eta+1}u_k
$$
converge absolutely in $L^p(M)$. Recall that $c_k\sim_{+\infty}
k^{\eta-\frac{1}{2}}$. Consequently, the series
$$
\sum_{k=1}^{+\infty}c_kT^{k-1}(I-T^2)^{\eta}(I+T)^{\alpha}(I-T)^{}u_k
$$
also converge absolutely. Now, for any $y\in L^{p^*}(M)$, we have
\begin{align*}
\MoveEqLeft
\Big\langle(N-1)!(I-T)x_1,y\Big\rangle_{L^{p}(M),L^{p^*}(M)}
    = \Big\langle(N-1)!(I-T)^{}Z^*u,y\Big\rangle_{L^{p}(M),L^{p^*}(M)}\\
   &= \Big\langle u, (N-1)!Z(I-T^{*})^{}y \Big\rangle_{L^{p}(M,\ell^2_\rad),L^{p^*}(M,\ell^2_\rad)} \\
   &= \bigg\langle u,\Big(c_k(T^*)^{k-1}\big(I-(T^*)^2\big)^{\eta}(I+T^*)^{\alpha}(I-T^*)^{}y\Big)_{k\geq 1}\bigg\rangle_{L^{p}(M,\ell^2_\rad),L^{p^*}(M,\ell^2_\rad)}\\
   &= \sum_{k=1}^{+\infty}\Big\langle u_k,c_k(T^*)^{k-1}\big(I-(T^*)^2\big)^{\eta}(I+T^*)^{\alpha}(I-T^*)^{}y\Big\rangle_{L^{p}(M),L^{p^*}(M)}\\
   &= \Bigg\langle\sum_{k=1}^{+\infty} c_kT^{k-1}\big(I-T^2\big)^{\eta}(I+T)^{\alpha}(I-T)^{}u_k ,y\Bigg\rangle_{L^{p}(M),L^{p^*}(M)}.
\end{align*}
Thus, we deduce that
\begin{equation}\label{(I-T)x1}
(N-1)!(I-T)^{}x_1=\sum_{k=1}^{+\infty}
c_kT^{k-1}\big(I-T^2\big)^{\eta}(I+T)^{\alpha}(I-T)^{}u_k.
\end{equation}
We have a similar result for $x_2$.
Now, we infer that
\begin{align*}
\MoveEqLeft (N-1)!(I-T)^{}(x_1+x_2)
    = (N-1)!(I-T)^{}x_1+(N-1)!(I-T)^{}x_2  \\
   &= \sum_{k=1}^{+\infty}  c_kT^{k-1}\big(I-T^2\big)^{\eta}(I+T)^{\alpha}(I-T)^{}u_k +\sum_{k=1}^{+\infty} c_kT^{k-1}\big(I-T^2\big)^{\eta}(I+T)^{\alpha}(I-T)^{}v_k \\
   &= \sum_{k=1}^{+\infty}  c_kT^{k-1}\big(I-T^2\big)^{\eta}(I+T)^{\alpha}(I-T)^{}(u_k+v_k)\\
   &= \sum_{k=1}^{+\infty}  c_kT^{k-1}\big(I-T^2\big)^{\eta}(I+T)^{\alpha}(I-T)^{}k^{\alpha-\frac{1}{2}}T^{k-1}(I-T)^{\alpha}x\hspace{1cm} \text{by (\ref{uk+vk})}\\
   &= \sum_{k=1}^{+\infty}  c_kk^{\alpha-\frac{1}{2}}T^{2k-2}\big(I-T^2\big)^{\eta+\alpha}(I-T)^{}x.
\end{align*}
For any $z\in \mathbb{D}$, we have
$$
\sum_{k=1}^{+\infty}k(k+1)\cdots
(k+N-2)z^{k-1}=\frac{(N-1)!}{(1-z)^N}.
$$
As in the proof of \cite[Theorem 3.3]{ALM}, we deduce that for every
$\varrho\in ]0,1[$ we have
\begin{equation}\label{formule pour Id}
(N-1)!I=\sum_{k=1}^{+\infty}c_kk^{\alpha-\frac{1}{2}}(\varrho
T)^{2k-2}\big(I-(\varrho T)^2\big)^{\eta+\alpha},
\end{equation}
the series being absolutely convergent. Hence, for any $\varrho\in
]0,1[$, we have
\begin{align*}
(N-1)!(I-\varrho T)x
   &=(I-\varrho T)\sum_{k=1}^{+\infty}c_kk^{\alpha-\frac{1}{2}}(\varrho T)^{2k-2}\big(I-(\varrho T)^2\big)^{\eta+\alpha}x \\
   &= \sum_{k=1}^{+\infty}c_kk^{\alpha-\frac{1}{2}}(\varrho T)^{2k-2}\big(I-(\varrho T)^{2}\big)^{\eta+\alpha}(I-\varrho T)^{}x.
\end{align*}
It is not difficult to see that the latter series is normally
convergent. Hence, letting $\varrho$ to $1$, we deduce that
$$
(N-1)!(I-T)x=
\sum_{k=1}^{+\infty}c_kk^{\alpha-\frac{1}{2}}T^{2k-2}\big(I-T^2\big)^{\eta+\alpha}(I-T)^{}x.
$$
Then we obtain
$$
(I-T)(x_1+x_2)=(I-T)x.
$$
Since the range $R(I-T)$ is dense in $L^p(M)$, the operator $I-T$ is
injective. Consequently, we have $x=x_1+x_2$. Now it remains to
estimate $\norm{x_1}_{T,\alpha,r}$ and $\norm{x_2}_{T,\alpha,c}$.
According to (\ref{(I-T)x1}), we have
\begin{align*}
m^{\alpha-\frac{1}{2}}T^{m-1}(I-T)^{\alpha} (N-1)!(I-T)^{}x_1 &=
\sum_{k=1}^{+\infty} c_k m^{\alpha-\frac{1}{2}}T^{k+m-2}
\big(I-T^2\big)^{\eta+\alpha}(I-T)u_k.
\end{align*}
Recall that $c_k\sim_{+\infty} k^{\eta-\frac{1}{2}}$ and that
$\sum_{k=1}^{+\infty} \norm{u_k}_{L^p(M)}^2 < \infty$. Then, it is
easy to see that, for any integer $m$, the series
$$
y_{m}=\sum_{k=1}^{+\infty} c_km^{\alpha-\frac{1}{2}}T^{k+m-2}
(I-T)^{\alpha+\eta}u_k
 \ \ \ \ \ \text{and} \ \ \ \ \
\sum_{k=1}^{+\infty} c_km^{\alpha-\frac{1}{2}}T^{k+m-2}
\big(I-T^2\big)^{\alpha+\eta}u_k
$$
are absolutely convergent in $L^p(M)$ on $[0,1]$. Let $y\in
R(I-T^*)$. There exists $y'\in L^{p^*}(M)$ such that $y=(I-T^*)y'$.
Then, we have
\begin{align*}
\MoveEqLeft
\Big\langle m^{\alpha-\frac{1}{2}}T^{m-1}(I-T)^{\alpha}(N-1)!x_1, y \Big\rangle_{L^p(M),L^{p^*}(M)}\\
   &= \Big\langle m^{\alpha-\frac{1}{2}}T^{m-1}(I-T)^{\alpha}(N-1)!x_1,(I-T^*)y' \Big\rangle_{L^p(M),L^{p^*}(M)}\\
   &= \Big\langle m^{\alpha-\frac{1}{2}}T^{m-1}(I-T)^{\alpha}(N-1)!(I-T)x_1,y' \Big\rangle_{L^p(M),L^{p^*}(M)} \\
   &= \Bigg\langle\sum_{k=1}^{+\infty} c_k m^{\alpha-\frac{1}{2}}T^{k+m-2} (I-T^2)^{\eta+\alpha}(I-T)u_k, y'\Bigg\rangle_{L^p(M),L^{p^*}(M)}\hspace{1cm} \text{by (\ref{(I-T)x1})}\\
   &= \Bigg\langle\sum_{k=1}^{+\infty} c_k m^{\alpha-\frac{1}{2}}T^{k+m-2} (I-T^2)^{\eta+\alpha}u_k,(I-T^*)y'\Bigg\rangle_{L^p(M),L^{p^*}(M)}\\
   &= \Bigg\langle\sum_{k=1}^{+\infty} c_k m^{\alpha-\frac{1}{2}}T^{k+m-2} (I-T^2)^{\eta+\alpha}u_k,y\Bigg\rangle_{L^p(M),L^{p^*}(M)}.
\end{align*}
Since $R(I-T^*)$ is dense, by duality, we have
\begin{align}
\label{y_m} m^{\alpha-\frac{1}{2}}T^{m-1}(I-T)^{\alpha}(N-1)!x_1
    &= \sum_{k=1}^{+\infty} c_k m^{\alpha-\frac{1}{2}}T^{k+m-2}(I-T^2)^{\eta+\alpha}u_k\\
    &=(I+T)^{\eta+\alpha}y_{m}\nonumber.
\end{align}
Now, observe that
\begin{align*}
c_k m^{\alpha-\frac{1}{2}}T^{k+m-2} (I-T)^{\alpha+\eta} &=\frac{c_k
m^{\alpha-\frac{1}{2}}}{(k+m-1)^{\eta+\alpha}}\cdot
(k+m-1)^{\eta+\alpha}T^{k+m-2}\big(I-T\big)^{\eta+\alpha}.
\end{align*}
Recall that we have proved in the proof of \cite[Theorem 3.3]{ALM}
that the matrix
$$
\biggl[\frac{c_k
m^{\alpha-\frac{1}{2}}}{(k+m-1)^{\eta+\alpha}}\biggr]_{k,m\geq 1}
$$
represents an element of $B(\ell^2)$. Moreover, by Proposition
\ref{prop bornitude}, the set
$$
\Big\{ (k+m-1)^{\eta+\alpha}T^{k+m-2}(I-T)^{\eta+\alpha}\,:\,
k,m\geq 1 \Big\}
$$
is Col-bounded. By Proposition \ref{prop 2.5 version col}, we deduce
that $(y_m)_{m\geq 1} \in L^p\big(M,\ell^2_c\big)$ and that
\begin{align*}
\bnorm{(y_{m})_{m \geq 1}}_{L^p(M,\ell^2_c)}
&\lesssim \norm{u}_{L^p(M,\ell^2_c)}.
\end{align*}
It is crucial to note that in this estimate, the majorizing constant
hidden in the symbol $\lesssim$ does not depend of $u$ or $x$. Since
$\{T\}$ is Col-bounded, we have
\begin{align*}
\norm{x_1}_{p,T,c,\alpha}
   &= \Bnorm{\big(m^{\alpha-\frac{1}{2}}T^{m-1}(I-T)^{\alpha}x_1\big)_{m\geq 1}}_{L^p(M,\ell^2_c)}\\
   &= \frac{1}{(N-1)!}\Bnorm{\big((I+T)^{\alpha+\eta}y_m\big)_{m\geq1}}_{L^p(M,\ell^2_c)} \hspace{1cm} \text{by (\ref{y_m})} \\
   &\lesssim  \bnorm{(y_{m})_{m \geq 1}}_{L^p(M,\ell^2_c)}.
\end{align*}
Finally, we deduce that there exists a positive constant $C$ such
that
\begin{align*}
\norm{x_1}_{p,T,c,\alpha}  &\leq C \norm{u}_{L^p(M,\ell^2_c)}.
\end{align*}
We have a similar result for $x_2$. Finally, we have
\begin{align*}
\norm{x_1}_{p,T,c,\alpha}+\norm{x_2}_{p,T,r,\alpha}
    &\leq C\norm{u}_{L^p(M,\ell^2_c)}+C\norm{v}_{L^p(M,\ell^2_r)}\\
    &\leq C\norm{x}_{p,T,\alpha}.
\end{align*}
\end{proof}
PREUVE Avant de mettre alpha=eta=1

\begin{proof}
Let $x\in L^p(M)$ such that $\norm{x}_{p,T,\alpha}<\infty$. There
exists two elements $u\in L^p\big(M,\ell^2_c\big)$ and $v \in
L^p\big(M,\ell^2_r\big)$ such that for any positive integer $k$
\begin{equation}\label{uk+vk}
    u_k+v_k=k^{\alpha-\frac{1}{2}}T^{k-1}(I-T)^{\alpha}x
\end{equation}
and such that
$$
\norm{u}_{L^p(M,\ell^2_c)}+\norm{v}_{L^p(M,\ell^2_r)}\leq
2\norm{x}_{p,T,\alpha}.
$$
By Theorem \ref{thm principal ALM}, the square functions
$\norm{\cdot}_{p^*,T^*,\eta}$ and $\norm{\cdot}_{T^*,\beta}$ are
equivalent for any $\eta,\beta>0$. Then we can suppose $\eta$ such
that $\eta+\alpha$ is an integer $N\geq 1$.

Now, for any integer $k\geq 1$, we define the complex number
$$
c_k=\frac{k(k+1)\cdots (k+N-2)}{k^{\alpha-\frac{1}{2}}}
$$
with the convention that $c_k=\frac{1}{k^{\alpha-\frac{1}{2}}}$ if
$N=1$. We
have $c_k\sim_{+\infty} k^{\eta-\frac{1}{2}}$. 
It follows that, for any $y\in L^{p^*}(M)$, we have
\begin{align*}
\MoveEqLeft
\Bgnorm{\bigg(\frac{c_k(T^*)^{k-1}(I-T^*)^{\eta}(I+T^*)^{N}y}{(N-1)!}\bigg)_{k\geq 1}}_{L^{p^*}(M,\ell^2_{\rad})}\\
    &\lesssim \bgnorm{\frac{1}{(N-1)!}(I+T^*)^{N}}_{L^{p^*}(M)\to L^{p^*}(M)} \bgnorm{\Big(k^{\eta-\frac{1}{2}}(T^*)^{k-1}(I-T^*)^\eta y\Big)_{k\geq 1}}_{L^{p^*}(M,\ell^2_{\rad})}\\
    &\lesssim \norm{y}_{L^{p^*}(M)}\hspace{1cm} \text{by (\ref{dual square estimate})}.
\end{align*}
Then, we obtain the bounded linear map
$$
\begin{array}{cccc}
 Z:   &  L^{p^*}(M)   &  \longrightarrow   &  L^{p^*}\big(M,\ell^2_{\rad}\big)  \\
      &      y        &  \longmapsto       & \displaystyle \Bigg(\frac{c_k(T^*)^{k-1}(I-T^*)^{\eta}(I+T^*)^{N}y}{(N-1)!}\Bigg)_{k\geq 1}.  \\
\end{array}
$$
Recall that we have contractive inclusions
$L^{p}\big(M,\ell^2_{c}\big) \subset L^{p}\big(M,\ell^2_{\rad}\big)$
and $L^{p}\big(M,\ell^2_{r}\big) \subset
L^{p}\big(M,\ell^2_{\rad}\big)$. We define the elements $x_1$ and
$x_2$ of $L^p(M)$ by
$$
x_1=Z^*u \ \ \ \ \ \text{and} \ \ \ \ \ x_2=Z^*v.
$$
We will show that $x=x_1+x_2$. Since $T$ is a Ritt-operator, there
exists a positive constant $C$ such that
\begin{align*}
 \sum_{k=1}^{+\infty} \Bnorm{k^{\eta-\frac{1}{2}}T^{k-1}(I-T)^{\eta+1}}_{L^p(M) \to L^p(M)}^2
   &=  \sum_{k=1}^{+\infty}k^{2\eta-1}\bnorm{T^{k-1}(I-T)^{\eta+1}}_{L^p(M) \to L^p(M)}^2  \\
   &\leq  C^2\sum_{k=1}^{+\infty} k^{2\eta-1}\bigg(\frac{1}{k^{\eta+1}}\bigg)^{2}  \\
   &=C^2\sum_{k=1}^{+\infty}\frac{1}{k^{3}}< \infty.
\end{align*}
By \cite[Proposition 2.5]{JMX}, we have the contractive inclusion
$L^{p}\big(M,\ell^2_{c}\big) \subset \ell^2\big(L^p(M)\big)$. We
deduce that $\sum_{k=1}^{+\infty} \norm{u_k}_{L^p(M)}^2 < \infty.$
According to Cauchy-Schwarz inequality, we deduce that the series
$$
\sum_{k=1}^{+\infty}k^{\eta-\frac{1}{2}}T^{k-1}(I+T)^{N}(I-T)^{\eta+1}u_k
=
(I+T)^{N}\sum_{k=1}^{+\infty}k^{\eta-\frac{1}{2}}T^{k-1}(I-T)^{\eta+1}u_k
$$
converge absolutely in $L^p(M)$. Recall that $c_k\sim_{+\infty}
k^{\eta-\frac{1}{2}}$. Consequently, the series
$$
\sum_{k=1}^{+\infty}c_kT^{k-1} (I+T)^{N}(I-T)^{\eta+1}u_k
$$
also converge absolutely. Now, for any $y\in L^{p^*}(M)$, we have
\begin{align*}
\MoveEqLeft
\Big\langle(N-1)!(I-T)x_1,y\Big\rangle_{L^{p}(M),L^{p^*}(M)}
   = \Big\langle(N-1)!(I-T)^{}Z^*u,y\Big\rangle_{L^{p}(M),L^{p^*}(M)}\\
   &= \Big\langle u, (N-1)!Z(I-T^{*})^{}y \Big\rangle_{L^{p}(M,\ell^2_\rad),L^{p^*}(M,\ell^2_\rad)} \\
   &= \bigg\langle u,\Big(c_k(T^*)^{k-1} (I+T^*)^{N}(I-T^*)^{\eta+1}y\Big)_{k\geq 1}\bigg\rangle_{L^{p}(M,\ell^2_\rad),L^{p^*}(M,\ell^2_\rad)}\\
   &= \sum_{k=1}^{+\infty}\Big\langle  u_k,c_k(T^*)^{k-1} (I+T^*)^{N}(I-T^*)^{\eta+1}y\Big\rangle_{L^{p}(M),L^{p^*}(M)}\\
   &= \Bigg\langle\sum_{k=1}^{+\infty} c_kT^{k-1} (I+T)^{N}(I-T)^{\eta+1}u_k ,y\Bigg\rangle_{L^{p}(M),L^{p^*}(M)}.
\end{align*}
Thus, we deduce that
\begin{equation}\label{(I-T)x1}
(N-1)!(I-T)^{}x_1=\sum_{k=1}^{+\infty} c_kT^{k-1}
(I+T)^{N}(I-T)^{\eta+1}u_k.
\end{equation}
We have a similar result for $x_2$.
Now, we infer that
\begin{align*}
\MoveEqLeft (N-1)!(I-T)^{}(x_1+x_2)
   = (N-1)!(I-T)^{}x_1+(N-1)!(I-T)^{}x_2  \\
   &= \sum_{k=1}^{+\infty}  c_kT^{k-1}(I+T)^{N}(I-T)^{\eta+1}u_k +\sum_{k=1}^{+\infty} c_kT^{k-1}(I+T)^{N}(I-T)^{\eta+1}v_k \\
   &= \sum_{k=1}^{+\infty}  c_kT^{k-1}(I+T)^{N}(I-T)^{\eta+1}(u_k+v_k)\\
   &= \sum_{k=1}^{+\infty}  c_kT^{k-1}(I+T)^{N}(I-T)^{\eta+1}k^{\alpha-\frac{1}{2}}T^{k-1}(I-T)^{\alpha}x\hspace{1cm} \text{by (\ref{uk+vk})}\\
   &= \sum_{k=1}^{+\infty}  c_kk^{\alpha-\frac{1}{2}}T^{2k-2}(I+T)^{N}(I-T)^{N+1}x.
\end{align*}
For any $z\in \mathbb{D}$, we have
$$
\sum_{k=1}^{+\infty}k(k+1)\cdots
(k+N-2)z^{k-1}=\frac{(N-1)!}{(1-z)^N}.
$$
As in the proof of \cite[Theorem 3.3]{ALM}, we deduce that for every
$\varrho\in ]0,1[$ we have
\begin{equation}\label{formule pour Id}
(N-1)!I=\sum_{k=1}^{+\infty}c_kk^{\alpha-\frac{1}{2}}(\varrho
T)^{2k-2}\big(I-(\varrho T)^2\big)^{N},
\end{equation}
the series being absolutely convergent. Hence, for any $\varrho\in
]0,1[$, we have
\begin{align*}
(N-1)!(I-\varrho T)x
   &=(I-\varrho T)\sum_{k=1}^{+\infty}c_kk^{\alpha-\frac{1}{2}}(\varrho T)^{2k-2}\big(I-(\varrho T)^2\big)^{N}x \\
   &= \sum_{k=1}^{+\infty}c_kk^{\alpha-\frac{1}{2}}(\varrho T)^{2k-2}(I+\varrho T)^{N}(I-\varrho T)^{N+1}x.
\end{align*}
It is not difficult to see that the latter series is normally
convergent on [0,1]. Hence, letting $\varrho$ to $1$, we deduce that
$$
(N-1)!(I-T)x=
\sum_{k=1}^{+\infty}c_kk^{\alpha-\frac{1}{2}}T^{2k-2}(I+T)^{N}(I-T)^{N+1}x.
$$
Then we obtain
$$
(I-T)(x_1+x_2)=(I-T)x.
$$
Since the range $\Ran(I-T)$ is dense in $L^p(M)$, the operator $I-T$
is injective. Consequently, we have $x=x_1+x_2$. Now it remains to
estimate $\norm{x_1}_{T,\alpha,r}$ and $\norm{x_2}_{T,\alpha,c}$.
According to (\ref{(I-T)x1}), we have
\begin{align*}
 (N-1)!m^{\alpha-\frac{1}{2}}T^{m-1}(I-T)^{\alpha+1} x_1 &=
\sum_{k=1}^{+\infty} c_k m^{\alpha-\frac{1}{2}}T^{k+m-2}
(I+T)^{N}(I-T)^{N+1}u_k.
\end{align*}
Recall that $c_k\sim_{+\infty} k^{\eta-\frac{1}{2}}$ and that
$\sum_{k=1}^{+\infty} \norm{u_k}_{L^p(M)}^2 < \infty$. Then, it is
easy to see that, for any integer $m$, the series
$$
y_{m}=\sum_{k=1}^{+\infty} c_km^{\alpha-\frac{1}{2}}T^{k+m-2}
(I-T)^{N}u_k
 \ \ \ \ \ \text{and} \ \ \ \ \
\sum_{k=1}^{+\infty} c_km^{\alpha-\frac{1}{2}}T^{k+m-2}
\big(I-T^2\big)^{N}u_k
$$
are absolutely convergent in $L^p(M)$. Let $y\in R(I-T^*)$. There
exists $y'\in L^{p^*}(M)$ such that $y=(I-T^*)y'$. Then, we have
\begin{align*}
\MoveEqLeft
\Big\langle (N-1)!m^{\alpha-\frac{1}{2}}T^{m-1}(I-T)^{\alpha}x_1, y \Big\rangle_{L^p(M),L^{p^*}(M)}\\
   &= \Big\langle (N-1)!m^{\alpha-\frac{1}{2}}T^{m-1}(I-T)^{\alpha}x_1,(I-T^*)y' \Big\rangle_{L^p(M),L^{p^*}(M)}\\
   &= \Big\langle (N-1)!m^{\alpha-\frac{1}{2}}T^{m-1}(I-T)^{\alpha+1}x_1,y' \Big\rangle_{L^p(M),L^{p^*}(M)} \\
   &= \Bigg\langle\sum_{k=1}^{+\infty} c_k m^{\alpha-\frac{1}{2}}T^{k+m-2} (I+T)^{N}(I-T)^{N+1}u_k, y'\Bigg\rangle_{L^p(M),L^{p^*}(M)}\hspace{1cm} \text{by (\ref{(I-T)x1})}\\
   &= \Bigg\langle\sum_{k=1}^{+\infty} c_k m^{\alpha-\frac{1}{2}}T^{k+m-2} \big(I-T^2\big)^{N}u_k,(I-T^*)y'\Bigg\rangle_{L^p(M),L^{p^*}(M)}\\
   &= \Bigg\langle\sum_{k=1}^{+\infty} c_k m^{\alpha-\frac{1}{2}}T^{k+m-2} \big(I-T^2\big)^{N}u_k,y\Bigg\rangle_{L^p(M),L^{p^*}(M)}.
\end{align*}
Since $R(I-T^*)$ is dense, by duality, we have
\begin{align}
\label{y_m} (N-1)!m^{\alpha-\frac{1}{2}}T^{m-1}(I-T)^{\alpha}x_1
    &= \sum_{k=1}^{+\infty} c_k m^{\alpha-\frac{1}{2}}T^{k+m-2}(I-T^2)^{N}u_k\\
    &=(I+T)^{N}y_{m}\nonumber.
\end{align}
Now, observe that
\begin{align*}
c_k m^{\alpha-\frac{1}{2}}T^{k+m-2} (I-T)^{N} &=\frac{c_k
m^{\alpha-\frac{1}{2}}}{(k+m-1)^{N}}\cdot
(k+m-1)^{N}T^{k+m-2}\big(I-T\big)^{N}.
\end{align*}
Arguing as in the proof of \cite[Theorem 3.3]{ALM}, we see that the
matrix
$$
\biggl[\frac{c_k
m^{\alpha-\frac{1}{2}}}{(k+m-1)^{N}}\biggr]_{k,m\geq 1}
$$
represents an element of $B(\ell^2)$. Moreover, by Proposition
\ref{prop bornitude}, the set
$$
\Big\{ (k+m-1)^{N}T^{k+m-2}(I-T)^{N}\,:\, k,m\geq 1 \Big\}
$$
is Col-bounded. By Proposition \ref{prop 2.5 version col}, we deduce
that $(y_m)_{m\geq 1} \in L^p\big(M,\ell^2_c\big)$ and that
\begin{align*}
\bnorm{(y_{m})_{m \geq 1}}_{L^p(M,\ell^2_c)}
&\lesssim \norm{u}_{L^p(M,\ell^2_c)}.
\end{align*}
It is crucial to note that in this estimate, the majorizing constant
hidden in the symbol $\lesssim$ does not depend of $u$ or $x$. Since
$\{T\}$ is Col-bounded, we have
\begin{align*}
\norm{x_1}_{p,T,c,\alpha}
   &= \Bnorm{\big(m^{\alpha-\frac{1}{2}}T^{m-1}(I-T)^{\alpha}x_1\big)_{m\geq 1}}_{L^p(M,\ell^2_c)}\\
   &= \frac{1}{(N-1)!}\Bnorm{\big((I+T)^{N}y_m\big)_{m\geq1}}_{L^p(M,\ell^2_c)} \hspace{1cm} \text{by (\ref{y_m})} \\
   &\lesssim  \bnorm{(y_{m})_{m \geq 1}}_{L^p(M,\ell^2_c)}.
\end{align*}
Finally, we deduce that there exists a positive constant $C$ such
that
\begin{align*}
\norm{x_1}_{p,T,c,\alpha}  &\leq C \norm{u}_{L^p(M,\ell^2_c)}.
\end{align*}
Moreover, we have a similar result for $x_2$. Finally, we have
\begin{align*}
\norm{x_1}_{p,T,c,\alpha}+\norm{x_2}_{p,T,r,\alpha}
    &\leq C\norm{u}_{L^p(M,\ell^2_c)}+C\norm{v}_{L^p(M,\ell^2_r)}\\
    &\leq C\norm{x}_{p,T,\alpha}.
\end{align*}
\end{proof}
